\overfullrule0pt
\font\small = cmr7

\newcount\footnum
\footnum=1
\def\newftnum{\global\advance\footnum by 1}
\def\ft#1{\footnote{$^{\vrule width .2cm\the\footnum\vrule width .2cm}$}{#1}
\newftnum}


\newcount\minutes
\newcount\scratch

\def\timestamp{%
\scratch=\time
\divide\scratch by 60
\edef\hours{\the\scratch}
\multiply\scratch by 60
\minutes=\time
\advance\minutes by -\scratch
\the\day$.\the \month.\the\year$ at $\,$\hours:\null
\ifnum\minutes< 10 0\fi
\the\minutes}
\mag=\magstep1
\hsize 5.4in
\vsize 7.6in

\def\qed{\hfill $ \sqcup\!\!\!\!\sqcap$\medskip}
\def\ni{\noindent}
\font\small = cmr7

\input epsf.sty


\medskip
\centerline{\bf A global theory of flexes of periodic functions}

\medskip
\vskip .5cm
\centerline{\it by Gudlaugur Thorbergsson and Masaaki Umehara}
\vskip 2cm
\footnote{}{\hskip -.5cm \small Mathematics Subject Classification (2000). Primary 51L15  Secondary 53C75, 53A15}

{ \ni{\bf Abstract.} {
For a real valued periodic smooth function $u$ on ${\bf R}$, $n\ge
0$, one defines the {\it osculating polynomial
$\varphi_s$} ({\it of order $2n+1$}) {\it at a point $s\in {\bf R}$ } to be
the
unique trigonometric polynomial of degree $n$,
whose value and first $2n$ derivatives at $s$ coincide with those of $u$ at
$s$.
We will say that a point $s$  is a {\it clean maximal flex} (resp.~{\it
clean 
minimal flex}) of the function $u$ on $S^1$ if and only if
$\varphi_s\ge u$ (resp.~$\varphi_s\le u$) and the preimage
$(\varphi-u)^{-1}(0)$ is connected.
We prove that {\it any smooth periodic function  $u$ has at least $n+1$
clean maximal flexes of order $2n+1$ and at least $n+1$ clean minimal flexes
of order $2n+1$.}
The assertion is clearly reminiscent of Morse theory and generalizes
the classical four vertex theorem for convex plane curves.
\par}
\vskip 1cm

\noindent
{\bf \S 1 \hskip 0.1in  Introduction} \medskip

For a real valued $C^{2n}$-function $u$ on $S^1={\bf R}/2\pi {\bf Z}$, $n\ge
0$, one defines the
{\it osculating polynomial
$\varphi_s$} ({\it of order $2n+1$}) {\it at a point $s\in S^1$} to be the
unique trigonometric
polynomial of degree $n$,
$$
\varphi_s(t)=a_0+a_1\cos t+b_1\sin t+\cdots+a_n\cos nt+b_n\sin nt,$$ whose
value and first $2n$ derivatives at $s$ coincide with those of $u$ at $s$.
If $u$ is $C^{2n+1}$ and the value and the first $2n+1$ derivatives of
$u$ and $\varphi_s$ coincide in $s$, i.e., if $\varphi_s$ hyperosculates $u$
in $s$, then we call $s$ a
{\it flex of $u$} ({\it of order $2n+1$}).
Notice that the order $2n+1$ of the osculating polynomials and flexes in the
definition above
has been chosen such that it coincides with the dimension of the space
$A_{2n+1}$ 
of trigonometric polynomials of degree
$n$. 
Notice also that a flex of order one, i.e.~the case $n=0$, is
nothing but a critical point.
The existence of $2n+2$ flexes of order $2n+1$
for any $C^{2n+1}$-function $u$ on $S^1$ is an easy consequence of the
well-known fact that
a function has at least $2n+2$ zeros if its Fourier
coefficients 
$a_i$ and $b_i$ vanish for $i\le n$; see Appendix A for a proof.
Here we will prove the much more difficult result that
there are
$2n+2$ such flexes satisfying the {\it global property} that the osculating
polynomials $\varphi_s$
in the flexes 
support  $u$, i.e., either $\varphi_s\le u$ or $u\le \varphi_s$. More
precisely, we will say that a
point $s$  is a {\it clean maximal flex} (resp.~{\it clean minimal flex}) of
a $C^{2n}$-function
$u$ on $S^1$ if and only if
$\varphi_s\ge u$ (resp.~$\varphi_s\le u$) and the preimage
$(\varphi-u)^{-1}(0)$ is connected.
This terminology is compatible with our definition of a flex, since it is
easy to see   that a clean
maximal (or minimal) flex is a flex if $u$ is $C^{2n+1}$.
\medskip

Our main result is the following theorem.\medskip

\ni{\bf 1.1 Theorem.}
{\it Let $u$ be a real valued $C^{2n}$-function on $S^1$ where $n\ge 1$.
Then $u$ has at least $n+1$
clean maximal flexes of order $2n+1$ and at least $n+1$ clean minimal flexes
of
order $2n+1$.}\medskip

The theorem is not true if $n=0$. A continuous function $u$ on $S^1$ is
obviously
supported by constant functions in points where $u$ takes on its maximum and
minimum values, but it
does not have to be true that $u$ takes on its maximum and minimum value in
connected sets.

The above theorem is clearly reminiscent of Morse theory.  We would like to
point out a further
similarity. Assume that $u$ is a $C^2$-function on $S^1$ and define the
function $\varphi^\bullet_s$
for every $s\in S^1$ as the largest function in $A_3$ such that
$\varphi_s^\bullet\le u$ and
$\varphi_s^\bullet(s)=u(s)$. Typically, $u$ and $\varphi^\bullet_s$ have two
common values. A point
$s$ in $S^1$ is therefore exceptional if $u$ and $\varphi^\bullet_s$ have
only a common value in $s$
or if $u$ and $\varphi^\bullet_s$ have more than two common values. In the
first exceptional case we
have that $s$ is a minimal flex. We denote the number of such flexes (or
the corresponding
functions $\varphi_s^\bullet$) by
$s^\bullet$. Let
$t^\bullet$ denote the number of functions $\varphi_s^\bullet$ counted with
multiplicities having
more than two values in common with $u$. (If $\varphi_s^\bullet$ and $u$
have
$k$ values in common, the $\varphi_s^\bullet$ contributes $k-2$ to the
number $t^\bullet$.) If
$s^\bullet$ is finite, then $t^\bullet$ is finite too and the following
formula holds:
$$
s^\bullet-t^\bullet=2.
$$
A similar formula holds for the  functions  $\varphi^\circ_s$
 defined for every $s\in S^1$ as the smallest function in $A_3$ such that
$\varphi_s^\circ\ge u$
and $\varphi_s^\circ(s)=u(s)$. The two formulas taken together generalize
Theorem 1.1 for $n=1$. One
should expect such formulas to hold for every $n$, thus giving a
far-reaching generalization of
Theorem 1.1, but so far there is no such result. The above formula
implies a theorem on
strictly convex curves that was first proved by  Bose [Bo] in the generic
case and then
generalized by Haupt [Ha] to generic simple closed curves. It was proved for
general simple closed curves by the second author [Um] using
intrinsic circle systems, a method that will be generalized in
the present paper.

If $n=1$, the existence of four flexes of order three on a periodic function
does in fact
imply the so-called {\it four vertex theorem} for
strictly convex
curves in the
Euclidean plane. A smooth regular curve $\gamma$ has at any point $s$ an
{\it osculating circle} which
can be defined as the unique circle having at least a three point contact
with $\gamma$ in $s$. The
point
$s$ is called a {\it vertex of $\gamma$} if the osculating circle at $s$
has at least a four point
contact with $\gamma$ in
$s$, or, in other words, the osculating circle hyperosculates $\gamma$ in
$s$. The four vertex theorem
for strictly convex curves says that such curves have at least four
vertices. Theorem 1.1
now implies  that {\it there  are at least two vertices at
which the osculating circles are inscribed and at least two vertices at
which they are circumscribed.}
This result is  more generally true for any simple closed curve in the
Euclidean plane, see
[Kn], and also follows from the methods we use here, see [Um],
but in this generality the curves do not correspond to functions on $S^1$.

We now describe the connection between strictly convex curves and
periodic functions in more detail. Fix a point $o$ in the interior
 of a strictly convex curve $\gamma$ in the
$(x,y)$-plane.

\medskip
\centerline{\epsfxsize=2.6in \epsfbox{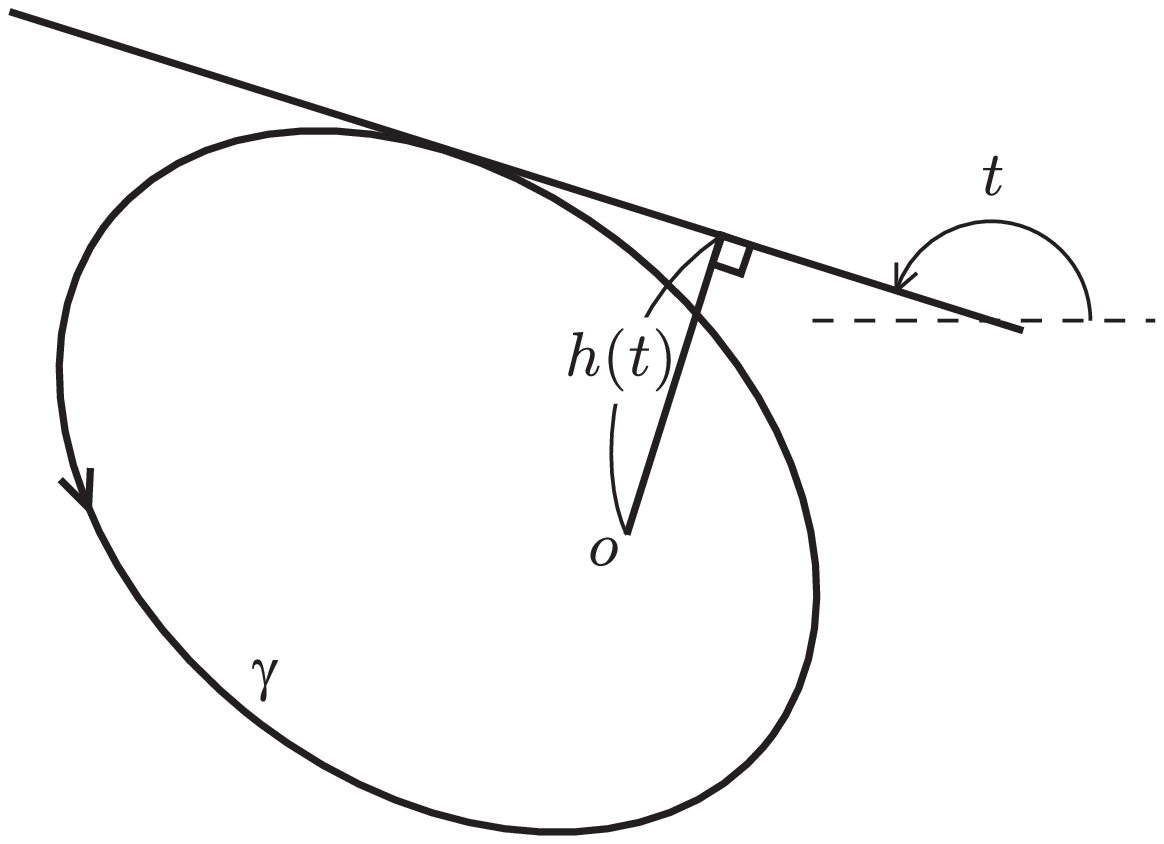}}
\centerline{Figure 1.}
\bigskip

\noindent For each $t\in [0,2\pi)$, there is a unique tangent
line $L(t)$ of the curve which makes angle $t$
with the $x$-axis. Let $h(t)$ be the distance between $o$ and the
line $L(t)$. The function $h$ is called the {\it supporting
function of the curve $\gamma$ with respect to $o$.}
The radius of the osculating circle of $\gamma$ at $t$ is given by
$h''(t)+h(t)$.
It can  easily be checked that a point $t$ is a vertex of the curve
if and only if $h'''(t)+h'(t)$ vanishes. Let $\varphi_s$ be the osculating
function (of order three) of $h$ in $s$. Then, by definition,
$h(s)=\varphi_s(s)$, $h'(s)=\varphi_s'(s)$ and $h''(s)=\varphi_s''(s)$.
Furthermore, $s$ is a flex if and only if $h'''(s)=\varphi_s'''(s)$.
Notice that $\varphi_s$ is a trigonometric polynomial of degree one, i.e.,
$$
\varphi_s(t)=a_0+a_1\cos t + b_1 \sin t.
$$
The function $\varphi_s'''+\varphi_s'$ clearly vanishes identically. Hence
$s$ is a flex of $h$ if and only if $h'''(s)+h'(s)=0$.
Hence we see that there is a one-to-one correspondence between
the vertices of the curve $\gamma$ and the flexes of the function $h$.
It is also easy to see that a clean maximal (resp.~minimal) flex
of $h$ corresponds to a vertex where the osculating circle is
inscribed (resp.~circumscribed ) and touches $\gamma$ in a connected set.

The methods of the paper are general and can be applied in other situations.
Let $\gamma$ again
denote a strictly convex closed curve that we assume to be contained in the
affine plane (or in the
projective plane, but this is not more general since there is always a line
which such a curve does
not meet). One defines the {\it osculating conic $C_t$ at a point $t$ of
$\gamma$} to be the unique
conic which meets $\gamma$ with multiplicity at least five in $t$. If $C_t$
and $\gamma$ meet with
multiplicity at least six in $t$, then $t$ is called a {\it sextactic
point.} If the osculating conic
at $t$ is inscribed (resp.~circumscribed) and meets $\gamma$ in a connected
set, then we will
call
$t$ a {\it clean maximal} (resp.~{\it clean minimal\/}) {\it sextactic
point}. One
can show that clean maximal and clean minimal sextactic points are in fact
sextactic.
We will prove the following theorem that improves a result of the authors in
 [TU2] as well as a result of Mukhopadhyaya in [Mu2].
\medskip

\ni{\bf 1.2 Theorem.}
{\it A strictly convex curve has at least three clean maximal sextactic
points
 and three clean
minimal sextactic points.}\medskip

The singularities we have been discussing so far are flexes of periodic
functions and vertices and sextactic points of convex curves. In all three
cases the dimension of the space
of approximating functions or curves is odd. (The space of circles is
three-dimensional and the space of
conics is five-dimensional.)  There are also singularities of even order of
which inflection points of
curves in the projective plane are the most important example. Here
one approximates the curve with
lines which form a two-dimensional space. In the even order case one
typically has to deal with
nonorientable situations like noncontractible curves in the projective plane
or antisymmetric functions. In
Appendix A we also deal with even order flexes of
antisymmetric functions, but we restrict
ourselves to the odd order case  in the main body of the paper since the
even order case is considerably
more difficult due to problems with nonorientability.

In Section two we introduce our main tool, the intrinsic systems,
which are in our applications analogues of intersection divisors on
algebraic curves.
 This approach is a
generalization of the methods in [Um], [TU1], [TU2], and
the idea behind it
is inspired by the paper [Kn] of Kneser. The
situation we deal with in Section two and in the rest of the paper is
somewhat more general than in the
introduction since we treat osculating functions that do not need to be
trigonometric functions.
In Section three we start drawing consequences from the defining axioms of
an intrinsic system.
In Section four we
generalize a result of Jackson [Ja].
In Section five we prove Theorem 1.1 and Theorem 1.2. Section six contains a
remark and some
questions on the
possible arrangement of the clean maximal and minimal flexes whose existence
was proved in Section
five. In Appendix A some basic properties of
trigonometric polynomials are
explained in the more general setting of Chebyshev spaces. In Appendix B we
explain an
elementary analytic result that is used in the paper.

\bigskip
\noindent
{\bf \S 2 \hskip 0.1in  Chebyshev spaces and intrinsic systems of
periodic functions }
\medskip

A real valued continuous function $u$ on $S^1$ is said to be
{\it piecewise $C^{2n}$} if it is of class $C^{2n}$ except at
finitely many  points $s_1,\dots,s_m$, and if,
furthermore, $u|_{[s_i,s_{i+1}]}$ can be extended to
a $C^{2n}$-function on
an open interval containing
$[s_i,s_{i+1}]$ for all $i=1,\dots,m$, where we understand
$m+1$ to mean $1$. We will refer to $s_1,\dots,s_m$ as {\it singular points}
or {\it
singularities of $u$}.

Our goal is to study the existence of clean flexes of order $2n+1$ of a
$C^{2n}$-function $u$ on $S^1$ that does not have
any singularities.
In the proofs below, we will frequently have to modify the function
$u$ by restricting it to an interval and then extending it to the
complement of the interval by piecewise trigonometric polynomials.
So we shall frequently have to deal with piecewise $C^{2n}$-periodic
functions. 
\medskip
We let $A_{2n+1}$ denote the vector space of
trigonometric polynomials of degree at most $n$, i.e.,
$$
A_{2n+1}=\left\{\varphi(t)= a_0+\sum_{k=1}^n (a_k\cos kt+b_k\sin kt)\right
\}. 
$$
The space $A_{2n+1}$ is an example of a  Chebyshev space of order $2n+1$;
see Appendix A where this concept is introduced and discussed in detail. We
repeat
here the definition of
Chebyshev spaces  of odd order.

\medskip
\noindent
{\bf 2.1 Definition.} A linear subspace ${\cal A}$ of $C^{2n}({\bf R}/2\pi)$
is called a {\it Chebyshev space of
order $2n+1$}  if its dimension is at least
$2n+1$ and if the number of zeros  in $[0,2\pi)$, counted with
multiplicities, of a nonvanishing function in ${\cal A}$ is at
most $2n$.

\medskip
It will be proved in Appendix A that  the dimension of a Chebyshev space
is always equal to its order. Let $u:S^1\to {\bf R}$
be a $C^{2n}$-function. Then for each $s\in S^1$
there exists a unique function $\varphi_s\in {\cal A}$
whose value and first $2n$ derivatives at $s$ coincide with those
of $u$ in $s$. We refer to Theorem A.2 in Appendix A for a proof of the
existence of $\varphi_s$.
We call 
$\varphi_s$  the {\it ${\cal A}$-osculating function} of $u$ at $s$.
If both $u$ and  $\varphi_s$ are $C^{2n+1}$-functions and
the value and the first $2n+1$ derivatives of
$u$ and $\varphi_s$ coincide in $s$,
then we call $s$ an {\it ${\cal A}$-flex of $u$}.\medskip

 We will from now on  work with an arbitrary Chebyshev space $\cal A$ of
order $2n+1$. The reader may want to think of $\cal A$ as simply
 being $A_{2n+1}$. Notice though that the more general
point of view is quite useful even when one is primarily interested
in $A_{2n+1}$. For an example of this, see the space ${\cal A}_{\psi_1}$
that is used to prove Theorem 1.1 from the introduction in Section five.
\medskip

Throughout the paper we let $I$ either denote the whole
$S^1$ or a nonempty
proper closed interval $[a,b]$ on $S^1$. In both
cases we will refer to
$I$ as an interval.

\medskip
\noindent
{\bf 2.2 Definition.}
Let $u$ be a  piecewise $C^{2n}$-function.
Let $I=[a,b]$ be a proper closed interval
on $S^1$ and $(\iota_a,\iota_b)$ a pair
of nonnegative integers which are less than or equal to $\infty$.
Then $u$ is said to satisfy the {\it boundary regularity condition
$(\iota_a,\iota_b)$} on $I$,
\item{(1)}  if  $u$ is at least $C^{2\iota_a-1}$ in $a$,
but not $C^{2\iota_a+1}$, and at least
$C^{2\iota_b-1}$ in $b$, but not $C^{2\iota_b+1}$,
\item{(2)}
and if $u$ is not $C^{2\iota_a}$ in $a$ (resp.~not $C^{2\iota_b}$
in $b$), then  the $2\iota_a$-th (resp.~$2\iota_b$-th) derivative
of $u$ from the left at $a$ (resp.~right at $b$) is greater than
that from the right at $a$ (resp.~left at $b$).

\medskip 

Let $I$ be a proper closed interval. We let $I^n_{(\iota_a,\iota_b)}$ denote
the subset of  the Cartesian product
$I^n$ consisting of those elements $(p_1,\dots,p_n)$
of $I^n$ with at most $\iota_a$ components equal to the
endpoint $a$ and at most $\iota_b$ components  equal to the
endpoint $b$. For example,
$$
\eqalign{
I^2_{(0,0)}&=(a,b)\times (a,b), \cr
I^2_{(1,0)}&=\{(x,y)\in [a,b]\times [a,b]\,;\, (x,y)\ne(a,a),\;x\ne b,\;
y\ne b\},\cr
I^2_{(1,1)}&=\{(x,y)\in [a,b]\times [a,b]\,;\, (x,y)\ne (a,a),(b,b) \}, \cr
I^2_{(2,1)}&=\{(x,y)\in [a,b]\times [a,b]\,;\, (x,y)\ne (b,b) \}, \cr
I^2_{(2,2)}&=[a,b]\times [a,b].}
$$

We next prove the following lemma.\medskip

\ni{\bf 2.3 Lemma.} {\it
Let $u$ be a piecewise $C^{2n}$-function on $S^1$
and let $I$ be a nonempty closed interval of
$S^1$ that is either proper or the whole circle. We suppose
that $u$ is $C^{2n}$ on $I$ and satisfies the  boundary regularity
 condition $(\iota_a,\iota_b)$ if
$I$ is a proper interval.
If $I= S^1$ we assume that  $u$ is $C^{2n}$ on
the whole $S^1$.
For $(p_1,\dots,p_n)\in I^n_{(\iota_a,\iota_b)}$  $($or
 $(p_1,\dots,p_n)\in
I^n$ if $I=S^1)$ let
$\Lambda$ denote the one-dimensional affine space of
functions $\varphi\in {\cal A}$ such that
$$
\varphi^{(k)}(p_i)=u^{(k)}(p_i)\quad {\rm for\; all}
\quad k=0,\dots,2\mu_i-1\,\,{\rm and\; all}\;\;i=1,\dots,n,
$$
where $\mu_i$ is the number of components
of $(p_1,\dots,p_{n})$ equal to
$p_i$. Then  the subset of
functions $\varphi \in \Lambda$ such that
$\varphi\ge u$ is a nonempty closed interval that we denote by
$\Lambda_u(p_1,\dots,p_n)$. }\medskip

\ni{\it Proof.} 
By definition, $\varphi\in \Lambda$ if and only if
$$
\varphi^{(k)}(p_i)=u^{(k)}(p_i)\quad {\rm for\; all}
\quad k=0,\dots,2\mu_i-1\,\,{\rm and\; all}\;\;i=1,\dots,n.
$$
It follows from Theorem A.2 in Appendix A
that $\Lambda$ is a one-dimensional
affine subspace of $\cal A$.
Let $\varphi_1$ be an arbitrary function in $\Lambda$.
Take another function $\varphi_2\in {\cal A}$ satisfying
$$
\varphi_2^{(k)}(p_i)=0 \quad {\rm for\; all}
\quad k=0,\dots,2\mu_i-1\,\,{\rm and\; all}\;\; i=1,\dots,n
$$
which is not identically zero. Then $\varphi_2^{(2\mu_i)}(p_i)
\ne 0$ for all $i=1,\dots,n$ (cf.~Definition 2.1).
Notice that $\varphi_2$ cannot change sign in any of the points $p_i$
since its first nonvanishing derivative there is of an
even order.
Since the function $\varphi_2$ has at most $2n$ zeros counted with
multiplicities,
either $\varphi_2\ge 0$ or $\varphi_2\le 0$ holds.
So we may assume $\varphi_2\ge 0$.
Then $\varphi_2^{(2\mu_i)}(p_i)>0$ for all
$i=1,\dots,n$. For every natural number
$m \in {\bf N}$,  we define a function $v_m$
on $S^1$ by setting
$$
v_m(t)=- u(t)+\varphi_1(t) +m\,\varphi_2(t).
$$
There is an $m_0\in {\bf N}$ such that for all
$m\ge m_0$ we have that $v_m$ and its first $2\mu_i-1$
derivatives vanish in
$p_i$, but $v_m^{(2p_i)}(p_i)>0$, for all
$i=1, \dots, n$, except when $p_i$ is either
$a$ or $b$ where $u$ might only
be $C^{2\iota_a-1}$ or
$C^{2\iota_b-1}$ respectively.
In case $p_i$ is  $a$ (resp.~$b$) and $u$ is only
$C^{2\mu_i-1}$, i.e.,
$\mu_i=\iota_a$ (resp.~$\mu_i=\iota_b$), we choose $m_0$
sufficiently large so that the $2\mu_i$-th
derivative from the left
and from the right of
$v_m$ are both positive in $p_i$.
Hence there is a neighborhood of
$\{p_1,\dots,p_n\}$ on which
$v_m$ is nonnegative.
On the complement of this neighborhood, we have
that $\varphi_2$ is
bounded from below by a positive number.
Hence there is a $m_1\ge m_0$ such
that
 $v_m\ge 0$ for $m\ge m_1$. Therefore the function
$$
\varphi(t)=\varphi_1(t)+m\, \varphi_2(t)
$$
is in $\Lambda$ and satisfies $\varphi\ge u$. The rest of the
lemma is now
clear.\qed

\medskip

Let $u$ be a function as in Lemma 2.3.
Now we can begin to associate what will call an  {\it intrinsic system}
to the function $u$. This is easy if $I$ is equal to the whole circle $S^1$
and the definition
consists only of the two cases (i), (ii) and (iv) below. We will therefore
restrict ourselves to the more difficult case of
functions  that are $C^{2n}$ on a proper closed interval $I$ satisfying a
boundary regularity
condition
$(\iota_a,\iota_b)$. We have seen in Lemma 2.3
that the subset $\Lambda_u(p_1,\dots,p_n)$
of $\cal A$ is a nonempty closed interval for each
$(p_1,\dots,p_n)\in I^n_{(\iota_a,\iota_b)}$.
We define the function $\varphi_{(p_1,\dots,p_n)}\in
\Lambda_u(p_1,\dots,p_n)$
to be the boundary point of this interval, or, what is the same thing, as
$$
\varphi_{(p_1,\dots,p_n)}=\inf\{\psi
\in \Lambda_u(p_1,\dots,p_n)\}.
$$
We call $\varphi_{(p_1,\dots,p_n)}$ {\it the minimal
function of $u$  with respect to  $(p_1,\dots,p_n)$.}

We denote by ${\bf N}_0$ the set of nonnegative integers,
and denote by ${{\rm Map}}(S^1,2{\bf N}_0\cup \{\infty\})$
the set of maps from $S^1$ to $2{\bf N}_0\cup \{\infty\}$.
We define a map
$$
f_u:I^n_{(\iota_a,\iota_b)}\to {{\rm Map}}(S^1,2{\bf N}_0\cup \{\infty\}),
$$
by setting   
$$
f_u(p_1,\dots,p_n)(q)=0\leqno{({\rm i})}
$$
for any $q\in S^1$ such that $u(q)\ne
\varphi_{(p_1,\dots,p_n)}(q)$;
$$
f_u(p_1,\dots,p_n)(q)=2k\leqno({\rm ii})
$$
 if $q\in I^\circ$ (where $I^\circ$
denotes the interior of $I$), $u(q)=
  \varphi_{(p_1,\dots,p_n)}(q)$
and precisely $2k-1$ derivatives of $u$ and
$\varphi_{(p_1,\dots,p_n)}$
agree in $q$ and $k\le n$;
$$
f_u(p_1,\dots,p_n)(q)=2k\leqno({\rm iii})
$$
if $q=a$  (resp.~$b$),
$u(q)=\varphi_{(p_1,\dots,p_n)}(q)$
and the first $2k-1$  derivatives of
$u$ and $\varphi_{(p_1,\dots,p_n)}$ agree in $q$,
the $2k$-th derivative of
$\varphi_{(p_1,\dots,p_n)}$ is different
from the $2k$-th derivative of $u$
from the right 
in $a$ (resp.~the left in $b$) and
$k\le\iota_a\le n$
(resp.~$\le\iota_b\le n$);
$$
f_u(p_1,\dots,p_n)(q)=\infty\leqno({\rm iv})
$$
 if $q\in I^\circ$, $u(q)= \varphi_{(p_1,\dots,p_n)}(q)$ and
 more than
$2n-1$ derivatives of $u$ and $\varphi_{(p_1,\dots,p_n)}$
agree in $q$;
$$
f_u(p_1,\dots,p_n)(q)=\infty\leqno({\rm v})
$$
if $q=a$ (resp.~$q=b$) and $u(t)$ is $C^{2n}$ at
$q$, 
$u(q)=\varphi_{(p_1,\dots,p_n)}(q)$ and more than
 $2n-1$ derivatives
of $u$ and $\varphi_{(p_1,\dots,p_n)}$ agree in $q$;
$$
f_u(p_1,\dots,p_n)(q)=2\leqno({\rm vi})
$$
if $q\not\in I$ and  $u(q)=\varphi(p_1,\dots,p_n)(q)$.

This ends the definition of the map
$$
f_u:I^n_{(\iota_a,\iota_b)}
\to {\rm Map}(S^1,2{\bf N}_0\cup\{\infty\}).
$$\medskip

It will frequently be convenient to use the following notation:
$$
f_u(p^k,p_{k+1},\dots,p_n)=f_u(p,\dots,p,p_{k+1},\dots,p_n),
$$
$$
f_u(p^k,q^l,p_{k+l+1},\dots,p_n)=f_u(p,\dots,p,q,\dots,q,p_{k+l+1},\dots,p_n
),
$$
and so on.
We will denote the support of $f_u(p_1,\dots,p_n)$
by $F_u(p_1,\dots,p_n)$, i.e.,
$$
F_u(p_1,\dots,p_n)=\{r\in S^1\;|\; f_u(p_1,\dots,p_n)(r)>0\}.
$$
 The value of $f_u(p_1,\dots,p_n)$ at a point $r$ will be called
 the {\it multiplicity of $r$ with respect
to $f_u(p_1,\dots,p_n)$.} The sum over all values
of $f_u(p_1,\dots,p_n)$, which
can of course be infinite, will be called the {\it total multiplicity of
$f_u(p_1,\dots,p_n)$}.

A point $s$ in $S^1$ will be called
a {\it global ${\cal A}$-flex of $u$}
if its multiplicity with respect
to  $f_u(s^n)$ is $\infty$.
Notice that {\it a point $s\in I$ is a a global ${\cal A}$-flex
if and only if $\varphi_{(s^n)}$ is defined and
equal to the ${\cal A}$-osculating function $\varphi_s$
of $u$ at $s$.}
 In particular, a global ${\cal A}$-flex is
 an ${\cal A}$-flex when $u$ and $\varphi_s$ are both  $C^{2n+1}$.
However, the converse is not true.
In fact, it is clear that a global ${\cal A}$-flex $s\in I$
has the global property that
the osculating function $\varphi_s$ of $u$ at
$s$ is greater than or equal to $u$ over the whole circle $S^1$.
A global ${\cal A}$-flex $s$ is called a
{\it clean maximal flex} if the preimage $(\varphi_s-u)^{-1}(0)$
 is connected.
If $u$ is  $C^{2n}$ on $S^1$,
then we can also define the intrinsic system
$f_{(-u)}$.
A clean maximal ${\cal A}$-flex of $-u$
is called a {\it clean minimal flex}.
Phrased differently, {\it a point $s$ is a
clean maximal $($resp. minimal$)$
${\cal A}$-flex of $u$ if and only if
the osculating function $\varphi_s$ is greater $($resp.~less$)$
than or equal to $u$ and
the preimage $(\varphi_s-u)^{-1}(0)$ is connected.}
\medskip

\noindent{\bf Example.} We give here an example which shows the difference
between
${\cal A}$-flexes, global ${\cal A}$-flexes and clean ${\cal A}$-flexes when
${\cal A}=A_{2n+1}$.
Consider a $2\pi$-periodic smooth function $u(t)$ satisfying $0\le u(t)\le
1$ which is identically
$1$ on the closed interval $I=[2\pi/5,3\pi/5]$
and identically zero on the intervals $[0,\pi/5]$ and $[4\pi/5,2\pi]$.

\medskip
\centerline{\epsfxsize=2.2in \epsfbox{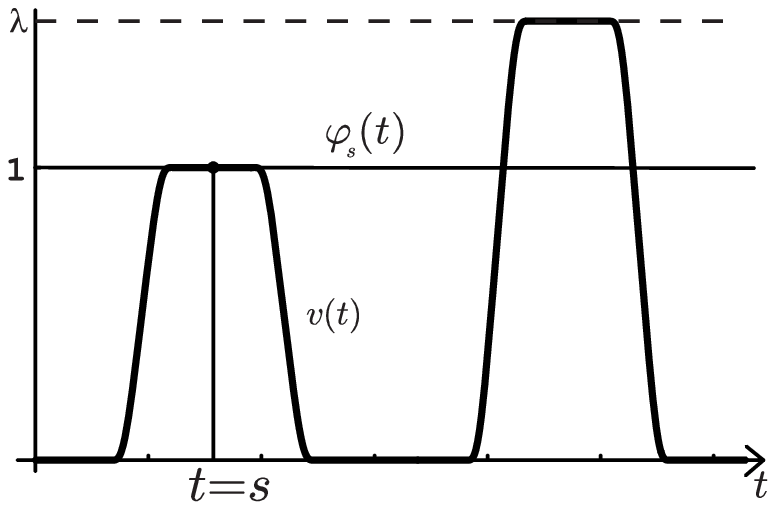}}
\centerline{A flex which is not clean}
\centerline{Figure 2.}

\bigskip
\noindent
Next we set
$$
v(t)=u(t)+ \lambda u(t+\pi)\qquad {\rm for}\quad \lambda\ge 0.
$$
When $\lambda<1$, the points on the interval $I$ are clean maximal
$A_{2n+1}$-flexes of $v(t)$.
If $\lambda=1$, the points on the interval $I$
are global $A_{2n+1}$-flexes
but not clean maximal $A_{2n+1}$-flexes.
Finally, if $\lambda>1$,
the points on the interval $I$
are ${A_{2n+1}}$-flexes, but not global $A_{2n+1}$-flexes.

\medskip

In the next proposition we bring the most basic properties of the map $f_u$
that we have associated to the function $u$. These properties will lead
us to the notion of an {\it intrinsic system}, see Definition 2.5 below.
Notice that
there are two cases. The interval $I$ is either a proper closed interval of
$S^1$ on which $u$ is
$C^{2n}$ and satisfies the boundary condition $(\iota_a,\iota_b)$ or $I$ is
the whole circle $S^1$
and
$u$ is $C^{2n}$ on $S^1$. We will formulate the following proposition for
the first case, i.e., for
$f_u:I^n_{(\iota_a,\iota_b)}
\to {\rm Map}(S^1,2{\bf N}_0\cup\{\infty\})$, but notice that everything is
equally true for
$I=S^1$; one simply has to delete the index $(\iota_a,\iota_b)$ from
$I^n_{(\iota_a,\iota_b)}$ and disregard (A.8).
\medskip

\noindent
{\bf 2.4 Proposition.} {\it
The map $f=f_u$ satisfies the following properties:
\medskip

\item{\rm (A1)} {\bf (Closedness)} The set $F(p_1,\dots,p_n)$ is closed
for all $(p_1,\dots,p_n)\in I^n_{(\iota_a,\iota_b)}$.

\item{\rm (A2)} {\bf (Symmetry)} The functions $f(p_1,\dots,p_n)$ and
$f(p_{\sigma(1)},\dots,p_{\sigma(n)})$
 coincide for every permutation
$\sigma\in {\cal S}_n$ and all
$(p_1,\dots,p_n)\in I^n_{(\iota_a,\iota_b)}$.

\item{\rm (A3)} {\bf (Supporting Property)} The support $F(p_1,\dots,p_n)$
contains $\{p_1,\dots,p_n\}$
for all $(p_1,\dots,p_n)\in I^n_{(\iota_a,\iota_b)}$.

\item{\rm (A4)} {\bf (Exchangeability)} If $f(p_1,\dots,p_n)(r)\ge 2j$ for a
point
$r\in I$ that is different from $p_1,\dots,p_{n-j}$ where $j\le n$, then
$(p_1,\dots,p_{n-j},r^j)\in I^n_{(\iota_a,\iota_b)}$ and
$f(p_1,\dots,p_{n-j},r^j)=f(p_1,\dots,p_n)$. In particular, if $j=n$,
then
$f(r^n)=f(p_1,\dots,p_n)$.

\item{\rm (A5)} {\bf (Uniqueness)}
If $p\in F(p_1,\dots,p_n)$ and $q\in F(q_1,\dots,q_n)$
satisfy
$$
p_1\preceq q_1\preceq\dots \preceq p_i\preceq
 q_i\prec p\prec q\prec p_{i+1}\preceq q_{i+1}\preceq\dots\preceq
p_n\preceq q_n(\prec p_1)
$$
where $0\le i\le n$, then $f(p_1,\dots,p_n)=f(q_1,\dots,q_n)$ holds.

\item{\rm (A6)} {\bf (Total Multiplicity)} The total multiplicity of
$f(p_1,\dots,p_n)$ is greater or
equal to
$2n+2$ for all $(p_1,\dots,p_n)\in I^n_{(\iota_a,\iota_b)}$
satisfying $f(p_1,\dots,p_n)(a)<2\iota_a$ and $f(p_1,\dots,p_n)
(b)<2\iota_b$.

\item{\rm (A7)} {\bf (Semicontinuity)} Let $(p_{1,k},\dots, p_{n,k})$
be a sequence in $I^n_{(\iota_a,\iota_b)}$  that converges to the element
$(p_1,\dots,p_n)\in I^n_{(\iota_a,\iota_b)}$ where $p_1\in I^\circ$.
Assume
$$
f(p_{1,k},\dots,p_{n,k})(p_{1,k})\ge 2\ell
$$
for all $k$. Then
$$
f(p_1,\dots,p_n)(p_1) \ge 2\ell.
$$
Assume $n\ge 2$. If $p_1=p_2$ and $p_{1,k}\ne p_{2,k}$ for all $k$, then
$$
f(p_1,\dots,p_n)(p_1) \ge 2\ell+2.
$$
\item{\rm (A8)} {\bf (Boundary isolation)}
If $\iota_a>1$ $($resp.~$\iota_b>1)$
and $0<f(p_1,\dots,p_n)(a)< 2\iota_a$ $($resp.~$0<f(p_1,\dots,p_n)(b)
< 2\iota_b)$, then $a$ $($resp.~$b)$ is isolated in $F(p_1,\dots,p_n)$.}
 \bigskip

\ni{\it Proof.}
Axioms (A1), (A2) and (A3) are trivially true for $f$.

We now prove that $f$ satisfies the Exchangeability Axiom (A4).
So we assume for $r\in I$ that
$$
f(p_1,\dots,p_n)(r)\ge 2j
$$
and $r\ne p_1,\dots,p_{n-j}$.
It follows from the definition of $f_u$ that
$(p_1,\dots,p_{n-j},r^j)\in I^n_{(\iota_a,\iota_b)}$.
We need to prove that $\varphi_{(p_1,\dots,p_{n-j},r^j)}=
\varphi_{(p_1,\dots,p_n)}$. It is clear
that
$\varphi(p_1,\dots,p_n)$ lies in $\Lambda_u(p_1,\dots,p_{n-j},r^j)$.
Hence
$$
\varphi_{(p_1,\dots,p_n)}
\ge\varphi_{(p_1,\dots,p_{n-j},r^j)}\ge u.
$$
Since we can squeeze $\varphi_{(p_1,\dots,p_{n-j},r^j)}$
between 
$\varphi_{(p_1,\dots,p_n)}$ and $u$ we have that
$\varphi_{(p_1,\dots,p_{n-j},r^j)}$ lies
$\Lambda_u(p_1,\dots,p_n)$. Hence
$$\varphi_{(p_1,\dots,p_{n-j},r^j)}\ge\varphi_{(p_1,\dots,p_n)}
\ge u.$$
It follows that
$\varphi_{(p_1,\dots,p_{n-j},r^j)}=\varphi_{(p_1,\dots,p_n)}$ and
hence 
$$
f(p_1,\dots,p_{n-j},r^j)=f(p_1,\dots,p_n),
$$
and  Axiom (A4) follows.

To prove that the Uniqueness Axiom (A5) is satisfied, assume that $p\in
F(p_1,\dots,p_n)$ and $q\in
F(q_1,\dots,q_n)$ satisfy
$$
p_1\preceq q_1\preceq\dots \preceq p_i\preceq
 q_i\prec p\prec q\prec p_{i+1}\preceq q_{i+1}\preceq\dots\preceq
p_n\preceq q_n(\prec p_1)
$$
where $0\le i\le n$.
Then the function
$$
\varphi_{(p_1,\dots,p_n)}-\varphi_{(q_1,\dots,q_n)}
$$ has more than
$2n+1$ zeros counted with multiplicities, implying that it  vanishes
identically. Hence we have $f(p_1,\dots,p_n)=f(q_1,\dots,q_n)$,
proving the axiom.

To prove Axiom (A6), first notice that the total multiplicity of $f$
is by definition greater or equal equal to $2n$.
Assume that the total multiplicity of $f$ is equal to $2n$.
Since $f(p_1,...,p_n)(a)<2\iota_a$
and $f(p_1,...,p_n)(b)<2\iota_b$,
 $u$ is  $C^{2\iota_a-1}$ at $a$
(resp.~$C^{2\iota_b-1}$ at $b$) and $a$ (resp.~$b$) occurs
less than $\iota_a$ times (resp.~$\iota_b$ times) as a component
of $(p_1,\dots,p_n)$.
The function $\varphi_{(p_1,\dots,p_n)}-u\ge 0$ and has precisely $2n$ zeros
counted with multiplicities.
More precisely, the set of zeros of $\varphi_{(p_1,\dots,p_n)}-u\ge 0$
is $\{p_1,\dots,p_n\}$, the first $2\mu_i-1$ derivatives
vanish in $p_i$ for all $i=1,\dots,n$, and the $2\mu_i$-th
derivative is positive in $p_i$.
Let $\varphi_2$ be as in the proof of Lemma 2.3.
The $2\mu_i$-th derivative of $\varphi_2$
is positive in $p_i$ for all $i=1,\dots,n$.
It follows that there is a sufficiently large $m$ such that
$$
\varphi_{(p_1,\dots,p_n)}-u\ge{1\over m}\varphi_2\ge 0.
$$
Hence  
$$
{\varphi_{(p_1,\dots,p_n)}}-{1\over m}\varphi_2\in
\Lambda_u(p_1,\dots,p_n),
$$
contradicting the definition of $\varphi_{(p_1,\dots,p_n)}$.

We now prove the Semicontinuity Axiom (A7).
 Let $(p_{1,k},\dots, p_{n,k})$
be a sequence in $I^n_{(\iota_a,\iota_b)}$  that converges to the element
$(p_1,\dots,p_n)\in
I^n_{(\iota_a,\iota_b)}$ where
$p_1\in I^\circ$.
Assume that
$$
f(p_{1,k},\dots,p_{n,k})(p_{1,k})\ge 2\ell
$$
for all $k$. Then clearly
$$
f(p_1,\dots,p_n)(p_1) \ge 2\ell,
$$
since $\varphi_{(p_{1,k},\dots,p_{n,k})}$ converges to
$\varphi_{(p_1,\dots,p_n)}$ together with all its derivatives.
Now assume that $n\ge 2$, $p_1=p_2$ and $p_{1,k}\ne p_{2,k}$ for all $k$.
We need
to prove that
$$
f(p_1,\dots,p_n)(p_1) \ge 2\ell+2.
$$
We consider the sequence $(v_k)$ where
$v_k=\varphi_{(p_{1,k},\dots,p_{n,k})}-u$.
 Notice that the first derivative of $v_k$
 vanishes in $p_{1,k}$ and
$p_{2,k}$ for all $k$. Hence
there is for every $k$ a point $q_{2,k}$ between $p_{1,k}$ and $p_{2,k}$
such that the second
derivative of $v_k$ vanishes in $q_{2,k}$. If $\ell\ge 2$,
then the second derivative of
$v_k$ vanishes in $p_{1,k}$ and $q_{2,k}$ and there must be a point
$q_{3,k}$ between $p_{1,k}$ and $q_{2,k}$ in which
the third derivative of $v_k$ vanishes. We continue this
argument inductively and show that there is for every $k$ a
point $q_{2\ell,k}$ between $p_{1,k}$ and $q_{2\ell,k}$ such
that the $2\ell$-th derivative of $v_k$ vanishes.
The sequence $(v_k)$ converges with all its derivatives
to $v=\varphi_{(p_{1},\dots,p_{n})}-u$ and the sequence
$(q_{2\ell,k})$ converges to $p$. It follows that at least
the $2\ell$ first derivatives of $\varphi(p_{1},\dots,p_{n})$
and $u$ coincide in $p_1$.  Since $u\le\varphi_{(p_{1},\dots,p_{n})}$
also the $2\ell+1$-st derivatives
coincide in $p$. Hence $f(p_{1},\dots,p_{n})(p_1)\ge 2\ell+2$ which
finishes the proof of Axiom (A7).

Finally we prove (A8). Suppose $\iota_a>1$ and
$0<f(p_1,...,p_n)(a)<2\iota_a$.
Suppose also that $a$ is not isolated in $F(a_1,...,a_n)$.
There is a sequence $(q_n)$ in $F(a_1,...,a_n)$
such that $\lim_{n\to \infty} q_n=a$. Then we have
$$
\varphi_{(p_1,...,p_n)}(q_n)=u(q_n)
$$
for all $n$.
Since $u(t)$ is $C^{2\iota_a-1}$ at $a$ this implies that
$$
\varphi_{(p_1,...,p_n)}(a)=u(a)
$$
and
$$
\varphi_{(p_1,...,p_n)}^{(j)}(a)=u^{(j)}(a)
$$
for $j=1,2,\dots, 2\iota_a-1$. It follows that
$f(p_1,...,p_n)(a)=2\iota_a$, a contradiction.
Hence $a$ is isolated in $F(a_1,...,a_n)$.
\qed

We now give the following definition.

\medskip
\noindent
{\bf 2.5 Definition.} Let $I$ either be the whole circle $S^1$
or a proper closed interval on $S^1$. In the second case we assume
we have a pair of $(\iota_a,\iota_b)$
of nonnegative integers which are less than or equal to $\infty$.
A map
$$
f:I^n_{(\iota_a,\iota_b)}\to {\rm Map}(S^1,2{\bf N}_0\cup\{\infty\}).
$$
is called an {\it intrinsic system of order $2n+1$
on $I$ (satisfying the boundary regularity
condition $(\iota_a,\iota_b)$)} if it satisfies
the axioms (A1) to (A8) in Proposition 2.4.
(If $I$ is the
whole circle one should of course delete
everything referring to the boundary conditions in the axioms.)
A point $s\in S^1$ is called an {\it $f$-flex}
if $f(s^n)(s)\ge 2n+2$.
Moreover, if $F(s^n)$ is connected, it is called a
{\it clean $f$-flex}.

\medskip
The map $f_u$ as in Proposition 2.4 is of course
an example of an intrinsic system of order $2n+1$, and
an $f_u$-flex is nothing but a global ${\cal A}$-flex.

Notice that the values of $f(p_1,\dots,p_n)$ can be finite numbers
greater than $2n$ although this does not happen for $f_u$ by
definition. This will for example happen in the course of the
reduction procedure introduced in Lemma 3.7
that we will frequently apply in the paper.

\medskip
We next give two more examples of intrinsic systems that come from curve
theory.

\medskip
\noindent{\bf 2.6 Example.} (i) Let $\gamma:S^1\to {\bf R}^2$ be
a simple closed regular $C^2$-curve.
For an arbitrary circle $C$,  we associate a function $\mu_C(r)$
on $S^1$ that maps a point $r$ on $\gamma$ to the multiplicity
with which $C$ and $\gamma$ meet in $r$.
The function $\mu_C(r)$ takes values in $\{0,1,2,3\}$
since we are only assuming the curve to be $C^2$-regular.
The value of $\mu_C(r)$ is of course zero in points
in which $C$ and $\gamma$ do not meet.
We let $C^\bullet_{p}$ ($p\in S^1$)
(resp.~$C^\circ_{p}$) denote the uniquely defined maximal
inscribed (resp.~minimal circumscribed) circle that is tangent
to $\gamma$ in $p$.  We set
$$
f_1^\bullet(p,q)(r)=\cases{\mu_{C^\bullet_{p,q}}(r) &
if  $\mu_{C^\bullet_{p,q}}(r)\le 2$,\cr
               \infty  & if $\mu_{C^\bullet_{p,q}}(r)\ge 3$.
               \cr
}
$$
We define the map $f_1^\circ$ similarly.
One can easily verify that $f_1^\bullet$ and $f_1^\circ$ are
both intrinsic systems of order $3$.
(Notice that the dimension of the space of
circles in the Euclidean plane is three.)
A point $p$ is called
a {\it clean maximal $($resp.~minimal$)$ vertex} if
the osculating circle $C_p$ is inscribed
(resp.~circumscribed) and meets the curve in a connected set.
When $\gamma$ is $C^3$,
the critical points of the curvature function of the curve
are called vertices. For a $C^3$-regular curve, clean vertices
are vertices of the curve. However
the vertices of a curve  have a priori no such global
properties. The concept of an intrinsic system is designed
to find vertices with such global properties.
It is well known and can be proved with the methods of this paper
that there are at least four clean vertices on a curve $\gamma$ as above.
The notion of an \lq intrinsic circle
system\rq\ as a family of closed subset $(F_p)_{p\in S^1}$ in $S^1$
satisfying certain axioms, see Section six below, was
introduced in [Um].
 Several applications were given  in [Um] and [TU1].
The family of supports $(F_1^\bullet(p))_{p\in S^1}$ and
$(F_1^\circ(p))_{p\in S^1}$ of the intrinsic systems $f_1^\bullet$ and
$f_1^\circ$
introduced above satisfy the axioms
of an intrinsic circle system.
As a consequence, when the curve has finitely many
maximal and minimal vertices,
one can prove that it satisfies a Bose type
formula as mentioned in the introduction; see [Um].

\medskip
(ii)
Let $\gamma$ be a strictly convex $C^4$-curve
in the real projective plane $P^2$. We identify $\gamma$ with $S^1$.
Let $\Gamma$ be a nondegenerate conic in $P^2$.
Then we associate to $\Gamma$ a function $\mu_\Gamma(r)$
on $S^1$ that maps a point $r$ on $\gamma$ to the multiplicity
with which $\Gamma$ and $\gamma$ meet in $r$.
The function $\mu_\Gamma(r)$ takes values in $\{0,1,\dots,5\}$
since we are only assuming the curve to be $C^4$-regular.
The value of $\mu_\Gamma(r)$ is of course zero in points
in which $\Gamma$ and $\gamma$ do not meet.
Let $(p,q)\in S^2$. If $p\ne q$, we let $\Gamma^\bullet_{p,q}$
(resp.~$\Gamma^\circ_{p,q}$) denote the uniquely defined maximal
inscribed (resp.~minimal circumscribed) conic that is tangent
to $\gamma$ in $p$ and $q$. If $p=q$, we let $\Gamma^\bullet_{p,q}$
(resp.~$\Gamma^\circ_{p,q}$) denote the uniquely defined maximal
inscribed (resp.~minimal circumscribed)
conic that meets $\gamma$ with multiplicity at least
four in $p=q$. We set
$$
f_2^\bullet(p,q)(r)=\cases{\mu_{\Gamma^\bullet_{p,q}}(r) &
if  $\mu_{\Gamma^\bullet_{p,q}}(r)\le 4$,\cr
               \infty  & if $\mu_{\Gamma^\bullet_{p,q}}(r)\ge 5$.
               \cr
}
$$
We define the map $f_2^\circ$ similarly.
One can easily verify that $f_2^\bullet$ and $f_2^\circ$ are
both intrinsic systems of order $5$.
(Notice that the dimension of the space of conics in $P^2$ is five.)
If an osculating conic at a point $p$ is inscribed
(resp.~circumscribed) and meets
the curve $\gamma$ in a connected set,
 we call $p$ a {\it clean maximal sextactic point}
 (resp.~{\it clean minimal sextactic point}).
When $\gamma$ is $C^5$, a point where the osculating conic meets with
multiplicity greater than $5$ is called a {\it sextactic
point}. 
By (A6), the clean maximal (resp.~minimal) sextactic points
are sextactic points of $\gamma$ whenever the curve is $C^5$.
 Existence of six sextactic points
where the osculating conics are inscribed or circumscribed
was proved by Mukhopadhyaya [Mu2];
see also [TU2] for an alternative proof.
These sextactic points might however not be clean.
We will refine the methods of [TU2] and prove
in Theorem 5.3 below
the existence of at least six clean sextactic points
on the curve $\gamma$.

In [TU2] we introduced the concept of an intrinsic conic system to
prove the above mentioned theorem of Mukhopadhyaya and more generally
 to find sextactic points on  simple closed curves in
$P^2$. Intrinsic conic systems are very similar  to
intrinsic systems of order five.\medskip
\medskip

We now generalize the construction of an intrinsic system
of order $2n+1$ on $I$
associated to a function $u$ by taking  {\it base points}
$c_{1},\dots,c_{r}$ into account since that will be
needed in Section five. Let $\nu_1,...,\nu_r$ be
positive integers. 
We set $N=n+m$ where $m=\sum_{h=1}^r \nu_h$ and let ${\cal
A}$  denote a Chebyshev space of order $2N+1$.
Let $I$ be a closed interval not containing
the base points $c_1,\dots,c_r$.
We assume we have a function $u$ that is
piecewise $C^{2n}$.
We assume furthermore that $u$ is $C^{2n}$ on $I$
satisfying the  boundary regularity condition $(\iota_a,\iota_b)$;
see Definition 2.2. We now generalize Lemma 2.3 to
this new situation.

\medskip
\ni{\bf 2.7 Lemma.} {\it
Assume the function $u$ and the base points $c_1,\dots,c_r$ to be as
described before this lemma.
We let $\mu_i$ $($resp. $\nu_{h})$ denote the multiplicity
with which $p_i$ $($resp.~$c_h)$ occurs as a component of the
$n$-uple 
$(p_1,\dots,p_n)$ $($resp.~$m$-uple $(c_1,\dots,c_m)\,)$.
Suppose also
that $u(t)$ be at least $C^{2\nu_h}$
on some neighborhood of $c_{h}$
for all $h=1,\dots,m$.
For $(p_1,\dots,p_n)\in I^n_{(\iota_a,\iota_b)}$ we let
$\hat \Lambda$ denote the one-dimensional set
of functions $\varphi\in {\cal A}$ such that
$$
\eqalign{
\varphi^{(k)}(p_j)&=u^{(k)}(p_j)
\qquad {\rm for }\;\,k=0,\dots,2\mu_j-1\;\;{\rm and}\;\;j=1,\dots,n \cr
\varphi^{(\ell)}(c_h)&=u^{(\ell)}(c_h)
\qquad\, {\rm for}\;\;\ell=0,\dots,2\nu_h-1\;\;{\rm and}\;\; h=1,\dots,m.
}
$$
Then  the subset of functions $\varphi \in \hat\Lambda$
such that $\varphi\ge u$ is a
nonempty closed interval that we denote by
$\hat\Lambda_{u}(p_1,\dots,p_n)$. }\medskip

\ni{\it Proof.}
One can proceed exactly as in the proof of Lemma 2.3.
\qed

\medskip
For a point $(p_1,\dots,p_n)\in I^n_{(\iota_a,\iota_b)}$
we define the function $\hat \varphi_{(p_1,\dots,p_n)}
\in \hat \Lambda_{u}(p_1,\dots,p_n)$ by setting
$$
\hat \varphi_{(p_1,\dots,p_n)}(t)=\inf\{\varphi(t)
\in \hat \Lambda_{u}(p_1,\dots,p_n)\}.
$$
As above we define a map
$$
\hat f_{u}:I^n\to {{\rm Map}}(S^1,2{\bf N}_0\cup \{\infty\}),
$$
satisfying the same first five conditions (i) to (v)
and the following three new conditions  (vi$^\prime$), (vii)
and (viii) (with (vi$^\prime$) replacing the previous
condition (vi)): 
$$
\hat f_{u}(p_1,\dots,p_n)(q)=2  \leqno({\rm vi'})
$$
if $q\not\in I\cup \{c_{1},\dots,c_r\}$ and  $u(q)=
\hat\varphi_{(p_1,\dots,p_n)}(q)$;
$$
\hat f_{u}(p_1,\dots,p_n)(q)=0 \leqno({\rm vii})
$$
if $q=c_h$ for some $h=1,\dots,r$ and precisely
 $2\nu_h-1$ derivatives of $u$ and $\hat\varphi_{(p_1,\dots,p_n)}$
 agree in $q$; and finally
$$
\hat f_{u}(p_1,\dots,p_n)(q)=2 \leqno({\rm viii})
$$
if $q=c_h$ for some $h=1,\dots,r$ and more than
 $2\nu_h-1$ derivatives of $u$ and $\hat \varphi_{(p_1,\dots,p_n)}$
 agree in $q$.\medskip

\medskip
\noindent
{\bf 2.8 Proposition.} {\it
The map $\hat f_u$ is an intrinsic system of order $2n+1$.
We shall call it the intrinsic system of order $2n+1$
 with base points $c_1, \dots c_r$ associated to $u$.}
\medskip

\noindent
{\it Proof.}
The proposition can be proved by modifying the
proof of Proposition 2.4.
\qed

\bigskip
\noindent
{\bf \S 3 \hskip 0.1in First consequences of the axioms of an intrinsic
system} \medskip

In this section, we shall derive some first
consequences of the axioms of intrinsic systems.
It should be remarked that Lemmas 3.3 to 3.6 below
are rather easy to check if the intrinsic system $f=f_u$
comes from a periodic function $u$.
Still we shall prove them only using the axioms since they are also
important 
for our applications to sextactic points in Section five.

The following trivial lemma  will frequently be used, mostly without saying
so explicitly.\medskip

\ni{\bf 3.1. Lemma.} {\it Let $J$ be a closed subinterval of the  interval
$I$.
Let $f$ be an intrinsic system of order $2n+1$ on $I$ satisfying the
boundary regularity
condition $(\iota_a,\iota_b)$
if $I$ is not the whole circle. Then the restriction of
$f$ to
$J^n\cap I^n_{(\iota_a,\iota_b)}$ is
an intrinsic system of order $2n+1$ on $J$
satisfying the
boundary regularity condition $(n,\iota_b)$
if $a$ is not in $J$ and $b$ is in $J$, the condition
$(\iota_a,n)$ if $b$ is not in
$J$ and $a$ is in $J$, and $(n,n)$ if neither $a$ nor $b$ lies in $J$.}
\qed

The following lemma is an immediate consequence of Axiom (A6).\medskip

\ni{\bf 3.2. Lemma.} {\it If $p\in I^\circ$
and $F(p,\dots,p)$
consists only of the point $p$, then $p$ is an $f$-flex.} \qed

The Exchangeability Axiom (A4) immediately implies the following
lemma.\medskip

\ni{\bf 3.3. Lemma.} {\it If $f(p_1,\dots,p_n)(p)\ge 2n+2$ for a point $p\in
I$, then $p$ is an
$f$-flex.}\qed

The next lemma is an application of the semicontinuity Axiom (A7).\medskip

\ni{\bf 3.4. Lemma.} (The Multiplicity Lemma) {\it We have
$f(p^j,p_{j+1},\dots,p_n)(p)\ge 2j$
for every $p\in I^\circ$.}\medskip

\ni{\it Proof.}  Let $(p_{l,k})$ for $l=1,\dots,j$ be $j$ sequences in
$I$ that converge to $p$ and assume that
$p_{l,k}\ne p$ for all
$l$ and all $k$. Axioms (A3) and (A7) imply that
$$f(p^2,p_{3,k},\dots,p_{j,k},p_{j+1},\dots,p_n)(p)\ge 4$$
for every $k$. (Here we fixed $k$ in the third and later arguments and let
$k$ in the first two
arguments go to infinity when applying (A7).)
We can now use the Symmetry Axiom (A2) to bring
$p_{3,k}$ into the second slot and use (A7) again to prove
$$f(p^3,p_{4,k},\dots,p_{j,k},p_{j+1},\dots,p_n)(p)\ge 6$$
for all $k$. We continue this argument inductively until we have proved the
lemma. \qed

\ni{\bf 3.5. Lemma.} {\it If $r\in F(p_1,\dots,p_n)\cap I^\circ$ is not
isolated in
$F(p_1,\dots,p_n)\cap I^\circ$, then
$r$ is an $f$-flex.}\medskip

\ni{\it Proof.} We assume that $f(p_1,\dots,p_n)(r)$ is a finite number
$k$. Set
$p_{1,k}=r$.
 Let $(p_{2,k})$
be a sequence in
$F(p_1,\dots,p_n)$ of pairwise different points that are all different from
$r$ and converge to $r$. After possibly permuting and relabeling the points
$p_1,\dots,p_n$, we have by the Exchangeability Axiom (A4) that
$f(p_1,\dots,p_n)(r)=f(p_{1,k},p_{2,k},p_3,\dots,p_n)(p_{1,k})=k$. Now the
Semicontinuity
Axiom (A7) implies that $f(p_1,\dots,p_n)(r)>k$, a contradiction. Hence
$f(p_1,\dots,p_n)(r)=\infty.$ It now follows from
the Exchangeability Axiom (A4) that $f(r^n)(r)=\infty.$\qed

The next lemma is the starting point of the idea of an
intrinsic system and the main tool in the paper [Um].
Notice that the Semicontinuity Axiom (A7) is
not used in its proof. The idea behind
the lemma goes back to H. Kneser [Kn].
Therefore we would like to call it the Kneser Lemma although
it is strictly speaking not due to him.
\medskip

\ni{\bf 3.6.  Lemma.} (The Kneser Lemma) {\it Let
$f$ be an intrinsic system of order three on
$I=[a,b]$ satisfying the boundary regularity
condition $(\iota_a,\iota_b)$
with
$\iota_a,\iota_b\ge 1$.
 Suppose that $a,b\in F(a)$ and $F(a) \cap
(a,b)$ is empty. Then there exists a point $c\in (a,b)$ such that
$F(c)$ is connected and contained in $(a,b)$. In particular, $c$ is an
$f$-flex.}\medskip

\ni{\it Proof.}
Let $q$ be any point in the interval
$(a,b)$. Then the Uniqueness Axiom (A5)  implies that
$F(q)$ is contained in $[a,b]$ and the Exchangeability Axiom
(A4) implies that $F(q)$ cannot contain
$ a$
and $b$.
Hence $F(q)\subset (a,b)$.
Let $c_1$ be the midpoint of the interval $[a,b]$.  If $F(c_1)$ is connected
then the proof is
finished. If $F(c_1)$ is not connected, there are two different points
$a_1,b_1\in F(c_1)$ such that
$F(c_1)\cap (a_1,b_1)$ is empty. Notice that the length of $[a_1,b_1]$
is less than half the length
of $[a,b]$ and $F(q)\subset(a_1,b_1)$ for every $q\in(a_1,b_1)$. Let $c_2$
be the midpoint of $[a_1,b_1]$. Then $F(c_2)\subset (a_1,b_1)$. If
$F(c_2)$ is connected we have finished the proof. If not, we continue
inductively and find a nested sequence of intervals $[a_n,b_n]$ with
midpoints $c_{n+1}$ such that $a_n,b_n\in F(c_n)$ and
$F(c_n)\cap(a_n,b_n)$ is empty. Furthermore, the length of $[a_n,b_n]$
is less than $(1/2)^n$ the length of $[a,b]$. We have
$F(q)\subset [a_n,b_n]$ for all $q\in (a_n,b_n)$. We stop the induction
if we arrive at a connected set $F(c_{n+1})$.  Otherwise we observe that
 the sequence $(c_n)$ converges to a point $c$. Then $F(c)$
consists of $c$ only since $F(c)\subset (a_n,b_n)$ for all $n$.
It now follows from  Lemma 3.2 that $c$ is an $f$-flex.
\qed

The main strategy in finding an $f$-flex of an intrinsic system of order
$2n+1$ is to reduce the order inductively until we can apply the Kneser
Lemma.
 We now start explaining this procedure.\medskip

Assume $n\ge2$ and let $f$ be an intrinsic system of order $2n+1$ on
$I=[a,b]$ satisfying
the boundary regularity condition $(\iota_a,\iota_b)$.
We choose $r=a$ (or $r=b$) and
assume
that $\iota_a\ge 2$ (or $\iota_b\ge 2$). If $r=a$,
let $(p_1,\dots,p_{n-1})\in
I^{n-1}_{(\iota_a-1,\iota_b)}$  and let $q\in S^1$. We set
$$
f_r(p_1,\dots,p_{n-1})(q)=\cases{f(r,p_1,\dots,p_{n-1})(q) &  if  $q\not
=
r$,\cr
                 f(r,p_1,\dots,p_{n-1})(r)-2 &  if  $q=r$, \cr}
$$
where we of course use the convention that $\infty-2=\infty$. We define
$f_r$
analogously on $I^{n-1}_{(\iota_a,\iota_b-1)}$ if $r=b$. \medskip

\noindent
{\bf 3.7. Lemma.} {\it
Let $I=[a,b]$ be a closed interval on $S^1$ and
$f$ an intrinsic system of order $2n+1$ on $I$ for some  $n\ge 2$
satisfying the boundary regularity condition $(\iota_a,\iota_b)$.
Let $r$ be an endpoint of $I$ and assume that $\iota_r\ge 2$.
Then $f_r$ is an intrinsic system of order $2n-1$ on $I$ satisfying
the boundary regularity condition $(\iota_a-1,\iota_b)$ if
$r=a$ and $(\iota_a,\iota_b-1)$ if $r=b$. }
\medskip

\noindent{\bf Remark.} The restriction in the Semicontinuity Axiom (A7)
that $p_1$ be in the interior $I^\circ$ of $I$ comes from
the fact that otherwise we would not be able to prove that
$f_r$ satisfies that axiom.
\medskip

\noindent{\it Proof.} We assume throughout the proof that $r=a$. The case
$r=b$ is completely analogous.

First notice that $f(r,p_1,\dots,p_{n-1})(r)\ge 2$ for
all $(p_1,\dots,p_{n-1})\in I^{n-1}_{(\iota_a-1,\iota_b)}$
by Axiom (A3). It follows
that the values of $f_r(p_1,\dots,p_{n-1})$ are nonnegative.

To see that (A1) is satisfied for $f_r$ we remark that the
sets $F(r,p_1,\dots,p_{n-1})$ and $F_r(p_1,\dots,p_{n-1})$ are
equal and
hence both closed if
$f(r,p_1,\dots,p_{n-1})(r)>2$. We have that
$F_r(p_1,\dots,p_{n-1})=F(r,p_1,\dots,p_{n-1})-\{r\}$ if
$f(r,p_1,\dots,p_{n-1})(r)=2$. Then (A8) implies
that $r$ is isolated in $F(r,p_1,\dots,p_{n-1})$ since $\iota_r\ge 2$.
Hence $F_r(p_1,\dots,p_{n-1})$ is also closed in this case.

Axioms (A2) and (A3) for $f$ clearly imply Axioms (A2) and (A3) for $f_r$.

To prove (A4) for $f_r$, assume that $f_r(p_1,\dots,p_{n-1})(q)\ge 2j$ for
$q\in I$
and $q$ does not coincide with any of the $p_1,\dots,p_{n-j-1}$ where $j\le
n-1$.
First assume that $q\ne r$. The  we have
$
f(r,p_1,\dots,p_{n-1})(q)\ge 2j
$
and Axiom (A4) for $f$ implies that
$$
f(r,p_1,\dots,p_{n-j-1},q^j)=f(r,p_1,\dots,p_{n-1}).
$$
Hence
$$
f_r(p_1,\dots,p_{n-j-1},q^j)=f_r(p_1,\dots,p_{n-1}).
$$
If $q=r$ we have $f(p_1,\dots,p_{n-1},r)(q)\ge 2j+2$.
 By Axiom (A4) for $f$ this implies that
$$
f(p_1,\dots,p_{n-j-1},q^j,r)=f(p_1,\dots,p_{n-1},r)
$$
Hence we again have that
$f_r(p_1,\dots,p_{n-j-1},q^j)=f_r(p_1,\dots,p_{n-1})$.

Axiom (A5) for $f_r$ follows immediately from Axiom (A5).

Axiom (A6) for  $f_r$ follows easily from Axiom (A6) for $f$ since
the total multiplicity of $f_r(p_1,\dots,p_{n-1})$ is two less than
the one of $f(r,p_1,\dots,p_{n-1})$ if the latter number is finite.
If the total multiplicity of $f(r,p_1,\dots,p_{n-1})$ is infinite,
then the same is true for $f_r(p_1,\dots,p_{n-1})$.

We now prove (A7) for $f_r$.  Let $(p_{1,k},\dots, p_{n-1,k})$
be a sequence in $I^{n-1}_{(\iota_a-1,\iota_b)}$ that converges to
the element $(p_1,\dots,p_{n-1})$ where $p_1\in I^\circ$. Notice that
$p_1\ne r$. Assume
$$
f_r(p_{1,k},\dots,p_{n-1,k})(p_{1,k})\ge 2\ell.
$$
Then
$$
f(p_{1,k},\dots,p_{n-1,k},r)(p_{1,k})\ge 2\ell
$$
for $k$ large. Now (A6) for $f$ immediately implies that
$$
f_r(p_1,\dots,p_{n-1})(p_1) \ge 2\ell.
$$
Assume $n\ge 3$. If $p_1=p_2\ne r$ and $p_{1,k}\ne p_{2,k}$ for all $k$,
then Axiom (A7) implies
$$
f_r(p_1,\dots,p_{n-1})(p_1)=f(p_1,\dots,p_{n-1},r)(p_1)\ge 2\ell+2.
$$
This shows that $f_r$ satisfies Axiom (A7).

Axiom (A8) for  $f_r$ follows easily from Axiom (A8) for $f$.
\qed

\noindent
{\bf 3.8. Lemma.} {\it Let $I=[a,b]$ be an interval on $S^1$,
$f$  an intrinsic system of order $2n+1$
on $I$ satisfying the boundary regularity condition $(n,\iota_b)$
with $\iota_b\ge1$.
Assume
that $f(a^n)=f(b,a^{n-1})$ and $F(a^n)\cap (a,b)$ is empty.  Then there
exists an
$f$-flex in the open interval $(a,b)$.

Similarly, if $f$ satisfies the boundary regularity
condition $(\iota_a,n)$ with
$\iota_a\ge1$,
$f(b^n)=f(a,b^{n-1})$ and $F(b^n)\cap (a,b)$ is empty, then there exists an
$f$-flex
in the open interval $(a,b)$. }
\medskip

\ni{\it Proof.} We proof the lemma by induction on $n$. The lemma is true
for $n=1$ by the Kneser Lemma 3.6. Assume the lemma is true for $n-1\ge 1$.

Assume that $f(a^n)=f(b,a^{n-1})$ and $F(a^n)\cap (a,b)$ is
empty where $f$ is
an intrinsic system of order $2n+1$. Then $f_a$ is an intrinsic system of
order $2n-1$ by Lemma 3.7.
Notice that $f_a(a^{n-1})=f_a(b,a^{n-1})$ and $F_a(a^{n-1})\cap
(a,b)=\emptyset$. By the
induction hypothesis there is a point $c\in (a,b)$ that is an
$f_a$-flex with respect to $f_a$. This implies that
$f(a,c^{n-1})(c)\ge 2n$.  By Axiom (A4) for $f$ this implies
$f(c^{n})=f(a,c^{n-1})$.
We can assume that $c$ is isolated in $F(c^n)$ since $c$ is otherwise
an  $f$-flex by Lemma 3.5 and there would be nothing left to prove.
Let $d$ be the point in $[a,c)\cap F(c^n)$ closest to $c$.
We have $f(c^n)=f(d,c^{n-1})$ and $F(c^n)\cap(d,c)=\emptyset$.
Then we can again use the induction hypothesis and we find
an  $f_c$-flex $e$ of $f_c$ in the interval $(d,c)$. Set $J=[e,c]$.

Let ${\cal C}$ denote the set of $(\alpha,\beta)\in J\times J$
such that $\alpha<\beta$, $f(\alpha^n)=f(\beta,\alpha^{n-1})$
and $F(\alpha^n)\cap (\alpha,\beta)=\emptyset$.
By arguments as in the previous paragraph we see that
${\cal C}$ is nonempty.

We let $\delta_{\alpha,\beta}$ denote  the distance between $\alpha$ and
$\beta$. Let $\delta$
denote the infimum over
$\delta_{\alpha,\beta}$ for
$(\alpha,
\beta)\in {\cal C}$.

We consider a sequence $\{(\alpha_k,\beta_k)\}$ in ${\cal C}$ such that
$\delta_{\alpha_k,\beta_k}$ converges to $\delta$.
By going to subsequences if necessary, we may assume that
$$
\lim_{k\to\infty}\alpha_k=\alpha,\qquad
\lim_{k\to\infty}\beta_k= \beta.
$$
If $\alpha=\beta$, then it follows immediately from Axiom (A7) since
$\alpha\in J\subset
I^\circ$ that
$$
f(\alpha^n)(\alpha)\ge 2(n+1)
$$
and we have that $\alpha\in J\subset (a,b)$ is an  $f$-flex.

We can therefore assume that $\delta>0$.
By (A7) we have
$f(\beta,\alpha^{n-1})(\alpha)\ge 2n$ and hence
$f(\alpha^n)=f(\beta,\alpha^{n-1})$ by the Exchangeability Axiom (A4).
We can assume that $\alpha$ and $\beta$ are isolated in $F(\alpha^n)$
since otherwise we have an $f$-flex by Lemma 3.5.
Let $\beta'$ be the point in $F(\alpha^n)\cap
(\alpha,\beta]$ closest to $\alpha$.  We now argue as in the second
paragraph of the proof  and
find  points $\gamma,\delta\in (\alpha,\beta')$ such that $(\gamma,
\delta)\in {\cal C}$.
Clearly $\delta_{\gamma,\alpha}<\delta$, which is a contradiction.
This finishes the proof of the
claim in the first paragraph of the lemma. The proof of the claim in the
second paragraph is similar.\qed
\medskip

The following two propositions are the main technical result of this
section.
Notice that very similar ideas go at least back to Mukhopadhyaya ([Mu1],
Propositions I and II) and Haupt and K\"unneth [HK], p.~47.
The main difference between our approach and theirs is that ours is more
global in nature and therefore allows us to prove the existence of flexes
satisfying global properties like being  clean.  The name of the
propositions
is taken from the book [HK].

\medskip
\noindent{\bf 3.9. Proposition.} (The Contraction Lemma I)
{\it Let $I=[a,b]$ be an interval on $S^1$,
$f$  an intrinsic system of order $2n+1$
on $I$ satisfying the boundary regularity condition $(\iota_a,\iota_b)$
with $\iota_a+\iota_b> n$.
Let $p_1,\dots,p_{n+1}\in I$ be such that $(p_1,\dots,p_n)\in
I^n_{(\iota_a,\iota_b)}$ and
$$
f(p_1,\dots,p_n)(a)+ f(p_1,\dots,p_n)(b)\ge 2(n+1).
$$
Then there exists an $f$-flex in the open interval $(a,b)$.
}
\medskip

\noindent{\sl Proof.}
We shall prove the proposition by induction.
If $n=1$, it follows from the Kneser Lemma 3.6.
We now assume that $n\ge 2$. Then  $\iota_a + \iota_b > n$ implies
that $\iota_a\ge 2$ or $\iota_b\ge 2$.
We consider the case $\iota_a\ge 2$ (the case $\iota_b\ge 2$
being similar). By Lemma 3.7, $f_a$ is an intrinsic system of
order $2n-1$ on $I$ satisfying the boundary regularity
condition $(\iota_a-1,\iota_b)$. By the induction hypothesis,
we find an $f_a$-flex $s$ on $(a,b)$. Then by the definition of
$f_a(s^{n-1})$, we have
$f(a,s^{n-1})(s)\ge 2n$ and hence $f(s^n)=f(a,s^{n-1})$.
If $s$ is not isolated in $F(a,s^{n-1})$, then $s$ is an $f$-flex by Lemma
3.5.
We can therefore assume that $s$ is isolated in $F(a,s^{n-1})$
and let $c$ be the point closest to $s$
in $[a,s)\cap F(a,s^{n-1})$. We get
$$
f(s^n)=f(c,s^{n-1})\quad {\rm and}\quad F(s^n)\cap (c,s)=\emptyset.
$$ 
Now by Lemma 3.8, we find an $f$-flex on $(c,s)\subset (a,b)$.
\qed

\medskip
\noindent{\bf 3.10. Proposition.} (The Contraction Lemma II) {\it Let
$f$ be an intrinsic system of order $2n+1$ on $I=[a,b]$ satisfying.
Let $p_1,\dots,p_{n+1}\in I^\circ$ be such that
 $p_1\preceq\dots\preceq p_{n+1}$ and
$$\sum_{t\in
\{p_1,\dots,p_{n+1}\}} f(p_1,\dots,p_n)(t)\ge 2(n+1).$$
(Notice that repeated points in the sequence $p_1,\dots,p_{n+1}$ only enter
once into the sum.)
Then there is an  $f$-flex in the open interval between $p_1$ and $p_{n+1}$.
}

\medskip

\noindent
{\bf Remark.} We do not assume any boundary condition in the proposition.
This is
possible since the points $p_1,\dots,p_{n+1}$ are assumed to be interior
points
of $I$.\medskip

\noindent{\sl Proof.}
We may assume $n\ge 2$ since the case $n=1$ follows easily from the Kneser
Lemma. 
We may also assume that $[p_1,p_{n+1}]\cap
F(p_1,\dots,p_n)$ consists of isolated points, since otherwise
there is either an  $f$-flex in $(p_1,p_{n+1})$ or we can find a smaller
interval with $n+1$ points from $F(p_1,\dots,p_{n+1})$ whose intersection
with $F(p_1,\dots,p_{n+1})$ consists of isolated points.

Assume that $p_1$ occurs $j$ times in the sequence  $p_1,\dots,p_{n+1}$. We
can assume that
$F(p_1,\dots,p_n)\cap (p_1,p_{j+1})$ is empty. If $j=n$ the claim
follows from Lemma 3.8. We therefore assume that $j<n$.
Then we can consider the intrinsic
system $g=f_{p_{n+1},p_n,\dots p_{j+3},p_{j+2}}$ of order $2j+1$ restricted
to
$[p_1,p_{j+1}]$ that we obtain by iterating the definition before
Lemma 3.7. This intrinsic
system satisfies the conditions in Lemma 3.8.
There is therefore a  $g$-flex $p_1'$
 in the open interval
$(p_1,p_{j+1})$. This implies
that $f(p_1,\dots,p_n)(p_1')\ge 2(j+1)$. We can therefore replace
$p_1,\dots,p_{j+1}$ by
$p_1'$ repeated $j+1$ times in the sequence $p_1,\dots, p_{n+1}$. We can
continue this
argument inductively until we are in the situation that $p_1$ occurs $n$
times in the sequence
$p_1,\dots, p_{n+1}$ and we can use Lemma 3.8 to find the
$f$-flex whose existence is claimed.
\qed

We now apply the methods of this section to prove a rather weak existence
theorem for $f$-flexes of an intrinsic system defined on the whole circle.
\medskip

\ni{\bf 3.11. Corollary.} {\it Let $f$ be an intrinsic system of order
 $2n+1\ge 5$ on $S^1$. Then $f$ has at least three  $f$-flexes.}\medskip

\ni{\bf Remark.} This assertion is optimal for $n=2$ as can be seen by
either considering sextactic points, see [TU2], or periodic functions,
see the example after Theorem 5.1. in Section five. We do not know whether
it is optimal for
$n>2$, but find it unlikely.  In the special case of intrinsic systems of
order $2n+1$ coming from periodic functions we will prove in Section five
the existence of at least $n+1$ points that are $f$-flexes, which is
optimal.
\medskip

\ni{\it Proof.}
We first prove the existence of two  $f$-flexes.
Let $p$ be some point on $S^1$ that is not an  $f$-flex. If such $p$ does
not exist there is nothing to prove. Notice that
$p$ is isolated in $F(p^n)$. Let $p_1$ and $p_2$ be the next points
to $p$ in $F(p^n)$ on each side of $p$.
It could happen that $p_1$ and $p_2$ coincide. By Lemma 3.8 there is
 an  $f$-flex in the open interval $(p_1,p)$ and another one in
the open interval $(p,p_2)$.

Now we prove the existence of a third $f$-flex. Denote the
$f$-flexes we have found by $q_1$ and $q_2$.
We first consider the possibility that $F(q_1,q_2^{n-1})$ only
consists of $q_1$ and $q_2$. Then it follows from the Contraction Lemma 3.9
(or 3.10)
that there is an $f$-flex in the open interval between
$q_1$ and $q_2$ and another one in the open interval between $q_2$ and
$q_1$.
 If $F(q_1,q_2^{n-1})$ has a point $q_3$ that is different
from $q_1$ and $q_2$, then we can after
renaming $q_1$ and $q_2$ if necessary assume that
$q_3$ lies in the open interval between $q_1$ and $q_2$. The Contraction
Lemma 3.10
then
implies the existence of an  $f$-flex in the open interval
between $q_1$ and $q_2$. This finishes the proof of the corollary.
\qed

\medskip

\bigskip
\noindent
{\bf \S 4 \hskip 0.1in The Jackson Lemmas} \medskip

We next prove two theorems - one for functions, the other for curves -
which we
will call  Jackson Lemmas, since a similar result
for vertices on simple  closed arcs was
found and applied by Jackson in [Ja], although
the existence of vertices having inscribed or circumscribed
osculating circles were not discussed in [Ja].
These two result will be used in Section five to prove the two theorems
stated in the introduction.

A piecewise $C^{2n}$-function %
$u$ will be  said to  have a {\it downward} (resp.~{\it upward\/})
 {\it singularity} at a singular point $s$, if
the interior angle in $s$ of the region
above (resp.~below) $u$ is less than or equal to $\pi$.

 \medskip
\ni{\bf 4.1 Theorem.} (The Jackson Lemma for Flexes of Functions) {\it
Let ${\cal A}$ be a Chebyshev space of order $2n+1$.
Let $u$ be a piecewise $C^{2n}$-function with at most one singularity
which we then denote by $a$. Suppose $u\not\in {\cal A}$ and that $a$ is
not an upward singular point of $u$.
Then $u$ has at least one clean  maximal ${\cal A}$-flex
with the property that the osculating function there does not have the same
value as $u$ in $a$.}\medskip

We first prove the following weaker version of the Jackson Lemma.

\medskip
\ni{\bf 4.2 Lemma.} {\it Let ${\cal A}$ and $u$ be as above.
Then $u$ has at least one global ${\cal A}$-flex $s$ such
that the osculating function $\varphi_s$ does not have the same
value as $u$ in $a$ and $\varphi_s\ge u$.}\medskip

\ni{\bf Remark.} The flex whose existence is claimed in Lemma 4.2 is an
$f_u$-flex with respect to the intrinsic system $f_u$ that can be associated
to $u$ (on a sufficiently large closed interval $I$ not containing
$a$).
\medskip

\ni{\it Proof.} We may assume that $a=0$.
We fix $n$ mutually distinct points $p_1\prec \dots\prec p_{n}$
arbitrarily, 
but all different from $a$. Since
$u$ is not in ${\cal A}$, we can assume that none of the points
$p_1,\dots,p_n$ is a
flex. We would like to show that the points $p_1,\dots,p_n$ can be
chosen such that $a\not\in F_u(p_1,\dots,p_n)$.
Assume  $a\in F_u(p_1,\dots,p_{n})$.
We choose points $q_1,\dots,q_n$ as follows:
$$
q_1\in (0,p_1),\quad q_2\in (p_1,p_2),\quad \dots,\quad q_n\in
(p_{n-1},p_n).
$$
Since $p_1,\dots,p_n$ are not flexes, it follows that $F_u(p_1,\dots,p_n)$
does not contain any of the intervals $(0,p_1),\dots, (p_{n-1},p_n)$
and hence also that we can choose $q_1,\dots,q_n$ such that
they are not contained in $F_u(p_1,\dots,p_n)$. (See Figure 3.
We indicate a periodic function $f(t)$ by the curve 
$\exp(u(t)-f(t))(\cos t,\sin t)$ on ${\bf R}$.) In particular, the function $f(t)=u(t)$ is expressed by the unit circle, and $f(t)=\varphi_u(p_1,\dots,p_{n})$
is expressed as a closed curve inscribed in the 
circle.

\medskip
\centerline{\epsfxsize=2.8in \epsfbox{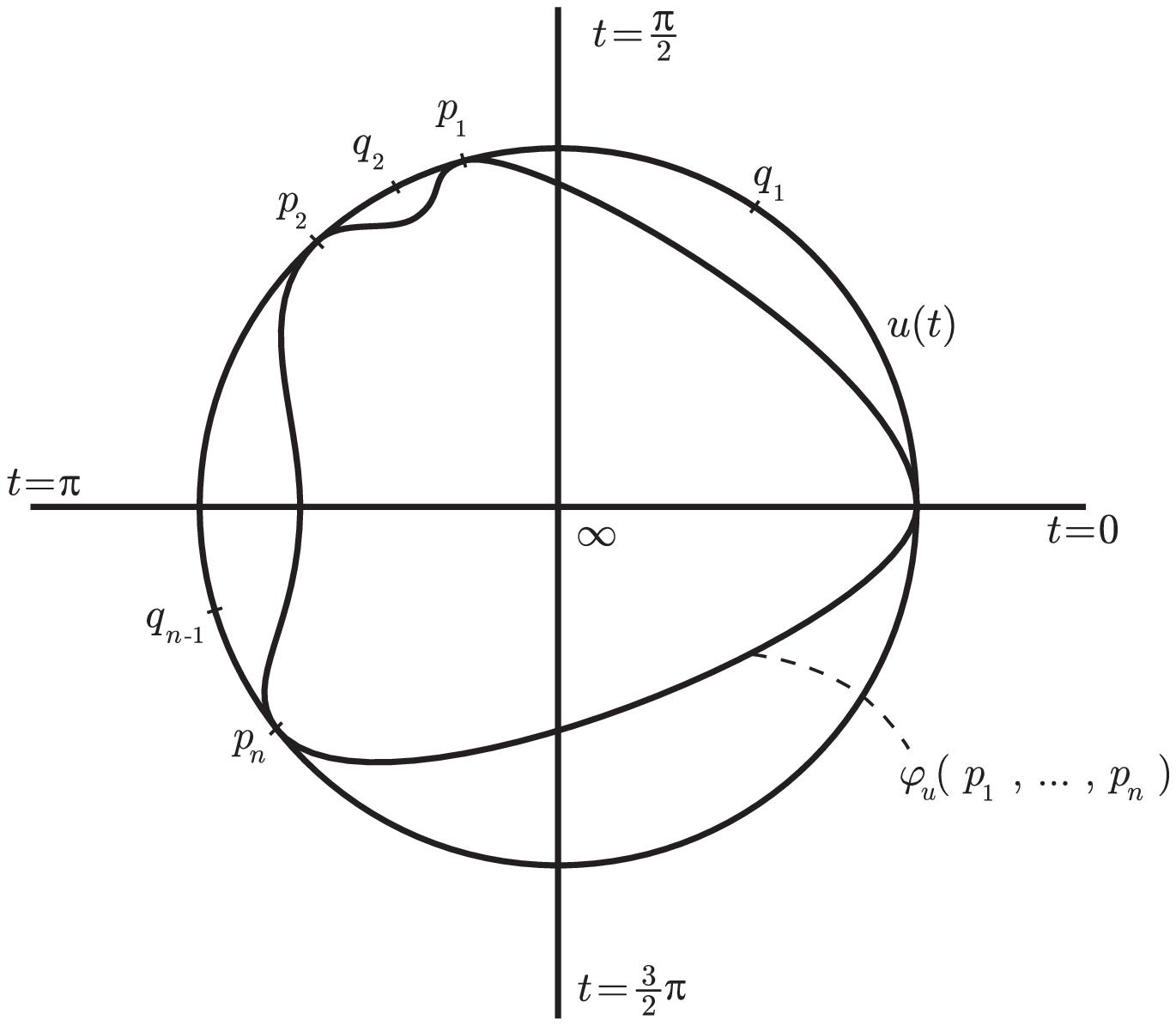}}
\centerline{Figure 3.}

\bigskip
\noindent
The graphs of $\varphi_u(p_1,\dots,p_{n})$ and $\varphi_u(q_1,\dots,q_{n})$
can at most meet in $2n$ points counted with multiplicities since
they are different. Hence they cannot meet in $a$. Now let $I$ be a closed
interval not containing $0$, but containing $F_u(p_1,\dots,p_n)$ in its
interior.
Since the total multiplicity
of $F_u(p_1,\dots,p_{n})$ is at least $2n+2$,
the Contraction Lemma 3.10 applied to $f_u$ restricted to $I^n$
implies the existence of a global ${\cal A}$-flex
$s$ on $I\subset (0,2\pi)$. Notice that we may assume that
the interval $[0,s]$ does not consist of ${\cal A}$-flexes.
We next show that there is such a flex with the property that
the osculating function does not have the same value as $u$ in $a$.

Suppose  that the osculating function $\varphi_{u}(s^n)$ has the same
value as $u$ in $a$. In this case, $u(t)$ must be $C^1$ at $0$.
Furthermore $\varphi_{u}(s^n)(t_0)
-u(t_0)>0$ holds  for some $t_0\in (0,s)$ since
the interval $[0,s]$ does not consist of ${\cal A}$-flexes.

We define
a function $v$ as follows:
$$
v(t)=\cases{u(t) & if $t\in [0,s]$, \cr
\varphi_u{(s^n)}(t)  & {\rm if $t\not\in [0,s]$}. \cr}
$$
Since $f_v$ is an intrinsic system of order $2n+1$ on $[0,s]$ satisfying
the boundary regularity condition $(1,n)$, we have
by Lemma 3.7 that $g:=(f_v)_{s^{n-1}}$ is an intrinsic system of
order $3$ on $[0,s]$ satisfying
$$
g(0)\ge 2, \quad g(s)\ge 2.
$$
Since $\varphi_{u}(s^n)(t_0)-u(t_0)>0$ holds  for some $t_0\in (0,s)$,
there is a point $r\in (0,s)$ such that $G(r)$ is connected
and $G(r)\subset (0,s)$ by the Kneser Lemma 3.6.
We define a new piecewise $C^{2n}$-function on $S^1$:
$$
w(t)=\cases{u(t) & if $t\in [r,s]$, \cr
\varphi_u{(r,s^{n-1})}(t)  & {\rm if $t\not\in [r,s]$}. \cr}
$$
Then we can define an intrinsic system $f_w$ of order $2n+1$
satisfying the boundary regularity condition $(2,n-1)$
on $[r,s]$.
Moreover we have
$$
f_w(r,s^{n-1})(r)\ge 4,\qquad f_w(r,s^{n-1})(r)\ge 2(n-1).
$$
Thus by the Contraction Lemma 3.9,
we find a global ${\cal A}$-flex $s'$ of $w$ on $(r,s)$, which
is a global ${\cal A}$-flex of $u$.
Since $G(r)\subset (0,s)$, we have $w(a)=\varphi_u{(r,s^{n-1})}(a)>u(a)$.
Thus the osculating function at $s'$ does not have the same
value as $u$ in $a$.
\qed

\noindent
{\it Proof of Theorem 4.1.}
We  let $\Phi(u)$ denote the set of global ${\cal A}$-flexes of $u$
with the property that the osculating functions at the flexes do not
have the same value as $u$ in $a$.
By Lemma 4.2, $\Phi(u)$ is nonempty.
For $p\in \Phi(u)$ we let $I(p)$ denote the minimal closed interval
in the complement of $a$ that contains $F_u(p^n)$.
We define a new piecewise $C^{2n}$-function $u_p$ on $S^1$
without upward singularities
by setting (See Figure 4.)
$$
u_p(t)=\cases{u(t) & if $t\in I(p)$, \cr
\varphi^u(p^n)(t) &  otherwise. \cr
}
$$

\smallskip
\centerline{\epsfxsize=2.8in \epsfbox{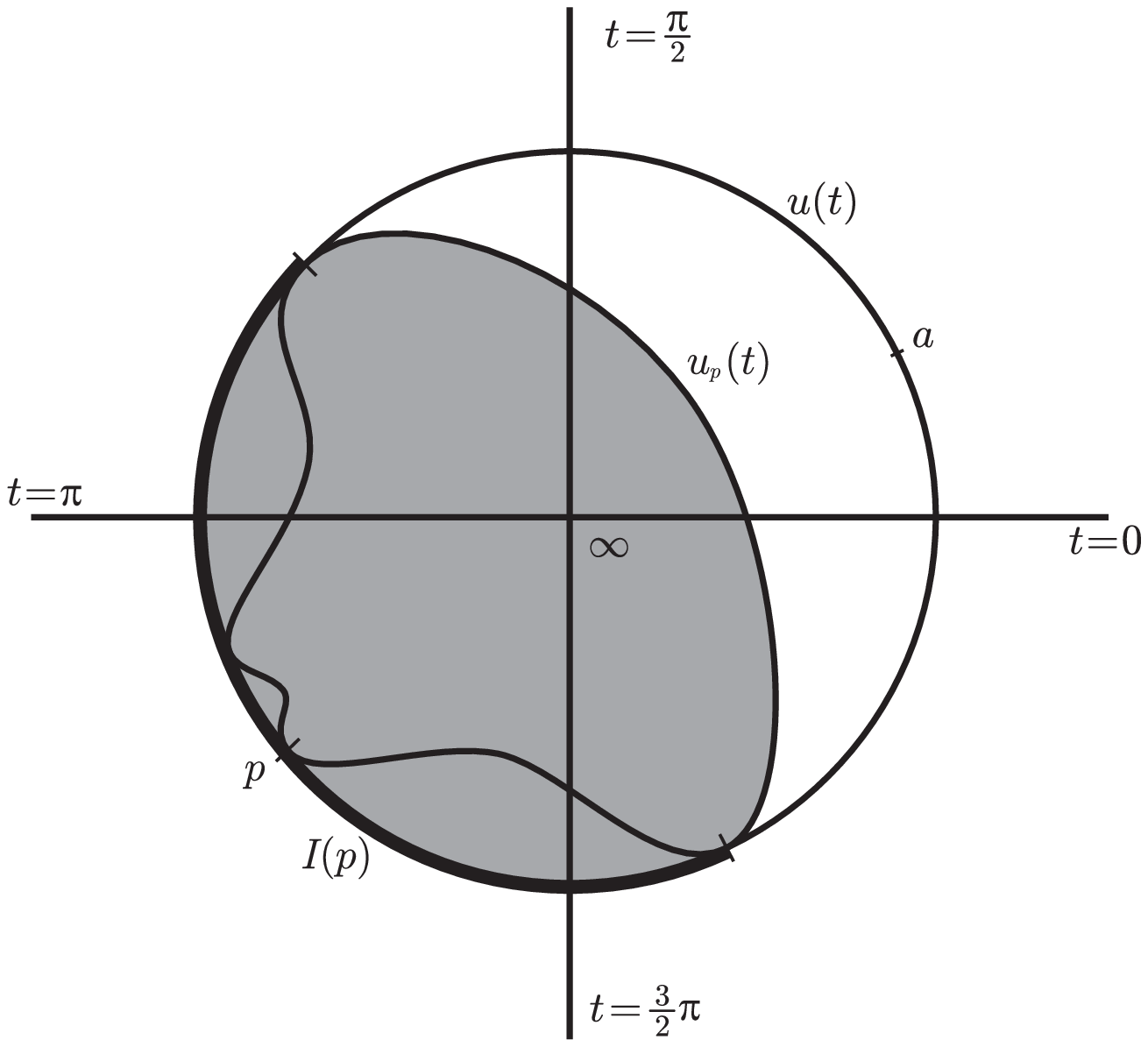}}
\centerline{Figure 4.}
\bigskip


We now define a partial ordering on $\Phi(u)$  by setting $p<\!\!<q$
for  $p,q\in \Phi(u)$ if
$$
I(p)\subset I(q) \qquad {\rm and } \qquad  u_q\le u_p.
$$
It is easy to check that this is in fact a partial ordering.

We next show that an element $p\in \Phi(u)$
that is minimal with respect
to this partial ordering is a clean flex, or,
in other words, $F_u(p^n)=I(p)$.
Assume that such a minimal element $p$ is not a clean flex of
$u$. Then $F_u(p^n)$ has at least two connected components.
Then there is a
point $q \in F_u(p^n)$ which belongs to a
connected component of $F_u(p^n)$  not containing $p$.
We can assume
that the open interval bounded by $q$ and $p$ does not
contain an $f_u$-flex. We
consider the case
$q\prec p$ (the case $p\prec q$ being similar).
We consider the piecewise $C^{2n}$-function $w$ on $S^1$
defined as
$$
w(t)=\cases{u(t) & if $t\in [q,p]$, \cr
\varphi_p(t)  & {\rm if $t\notin [q,p]$}, \cr}
$$
where $\varphi_p(t)$ is the osculating function of $u(t)$
at $p$.
The function $w$ is a $C^{2n}$-function on the interval $[q,p]$
satisfying the boundary regularity condition $(1,n)$,
and we see from the Contraction Lemma 3.9 that there exists
an $f_w$-flex $s$ in the open interval $(q,p)$.
The osculating function $\varphi_{w,s}$ of $w$ at $s$ is equal to
the osculating function $\varphi_{u,s}$ of $u$ at $s$ and
$$
\varphi_{u,s}=\varphi_{w,s}\ge w\ge u_p \ge u.
$$
Thus $s$ is also an $f_u$-flex. Notice that $s\in \Phi(u)$, since
$\varphi_{u,s}(a)\ge u_p(a)>u(a)$.
Furthermore $s<\!\!< p$ since
$\varphi_{u,s}(t)\ge \varphi_{u,p}(t)$ for
$t\not\in I(p)$. Since $s\ne p$,
this is a contradiction and we have proved that a  point $p$ which is
minimal with respect to the partial ordering is a clean flex.

We will now prove that there is a minimal point with respect to the partial
ordering with help of Zorn's Lemma.
Let $S$ be an arbitrary totally ordered subset of $\Phi(u)$.
We fix some point $p_0\in S$ and let $S_0$ denote
the set of elements $p\in S$ such that $p<\!\!< p_0$.
For $t\in S^1\setminus I(p_0)$, we set
$$
\varphi(t)=\sup\{\varphi_{u,p}(t)\;|\;p\in S_0\}.
$$
Notice that $\varphi_{u,p}$ depends continuously on $p$
since it is an osculating function
of $u$ and $u$ is $C^{2n}$.
Hence the family of osculating functions
of $u$ is bounded and the function $\varphi(t)$ is
well defined.
We would like to show that $\varphi(t)$ is the
restriction to $S^1\setminus I(p_0)$ of
a function in ${\cal A}$.
We fix $2n+1$ distinct points $t_1,\dots,t_{2n+1}$ on $S^1$
and set
$$
\alpha_j=\sup\{\varphi_{u,p}(t_j)\;|\;p\in S_0\}
\qquad {\rm for }\qquad j=1,\dots,2n+1.
$$
Then there exists a unique function
$\psi(t)\in {\cal A}$ such that
$$
\psi(t_j)=\alpha_j.
$$
There is a sequence $(p_k)$ in $S_0$ such that
$$
\alpha_j=\lim_{k\to\infty}\varphi_{u,p_k}(t_j)
\qquad {\rm for }\qquad j=1,\dots,2n+1.
$$
It follows that the sequence $(\varphi_{u,p_k})$ converges
uniformly to $\psi$ on $S^1$.
Suppose that $\psi(c)<\varphi(c)$ for some
$c\in S^1\setminus I(p_0)$.
Since $\varphi(c)=\sup\{\varphi_{u,p}(c)\,;\, p\in S_0\}$,
there exists $q\in S_0$ such that
$\psi(c)<\varphi_{u,q}(c)$.
In particular $\varphi_{u,p_k}(c)<\varphi(c)$.
Since $S_0$ is a totally ordered set, we have
$\varphi_{u,p_k}(t)\le\varphi_{u,q}(t)$ for all
$t\in S^1\setminus I(p_0)$. There is some $k_0$ such that
$q<p_{k_0}$. Hence $\varphi_{u,p_{k_0}}(t)\ge\varphi_{u,q}(t)$
and it follows that $\varphi_{u,p_k}(t)=\varphi_{u,q}(t)$ for $k\ge k_0$,
contradicting $\psi(c)<\varphi_{u,q}(c)$. It follows that $\varphi$ is a
restriction
of the function $\psi$ in $\cal A$ to $S^1\setminus I(p_0)$ as we wanted to
show.
We can assume that $(p_k)$ converges to a point $p_\infty$.
It is clear that
$p_\infty$ is a flex and that $\varphi(t)=\varphi_u(p_\infty)$.
Then it follows that
$q\in \Phi(u)$ and $p_\infty<\!\!<p$ for all $p\in S$
since $\varphi_u(p_\infty)$ is
strictly larger than $u$ outside of the interval $I$
defined by
$$
I=\cap_{p\in S} I(p).
$$
We can therefore use Zorn's Lemma to find a minimal point with respect to
the partial ordering
thereby proving the existence of the clean flex with the desired properties.
\qed

The following theorem is the analogue of the Jackson Lemma for sextactic
points. It will be used in section five to
prove the theorem on sextactic points from the introduction.

\medskip
\ni{\bf 4.3 Theorem.} (The Jackson Lemma for Sextactic Points) {\it
Let $\gamma:S^1\to {\bf R}^2$ be a simple closed curve which is not
a conic and is everywhere
$C^4$-regular except maybe in a given point $a$ where we assume that it
is $C^4$-regular from both left and
right.  We assume furthermore that $\gamma$ is bounds a convex
region and that it is strictly convex except maybe in the point $a$.
Then there is a clean sextactic  point $s$ on $(a,a+2\pi)$ with
the property that
the osculating conic at $s$ is inscribed
and does not meet $\gamma(a)$.  If furthermore,
$\gamma$ is at least $C^1$ in $a$, then $\gamma$ has a clean
sextactic point $s$ with the property that the osculating
conic in $s$ is circumscribed and does not meet $\gamma(a)$. }

\medskip
\ni{\it Proof.} In the proof we will assume that the curve $\gamma$ lies in
$P_2$, i.e., we compactify ${\bf R}^2$ by adding a line at infinity.
It was explained in Example 2.6 (ii) how a regular strictly convex curve
in  the affine plane gives rise to an intrinsic system.
Here the situation is somewhat different since we are
allowing a singular point $a$.

Without loss of generality, we may set $a=0$.
Take two distinct points $p_1,p_2\in (0,2\pi)$ which are not
sextactic points.
Consider a maximal inscribed conic $\Gamma^\bullet_{p_1,p_2}$
(resp.~a minimal circumscribed conic $\Gamma^\circ_{p_1,p_2}$)
passing through $p_1$ and $p_2$.
Suppose $\Gamma^\bullet_{p_1,p_2}$ passes through $a$. Then $\gamma$
must be $C^1$ at $a$.
Choose $q_1\in (0,p_1)$ and
$q_2\in (p_1,p_2)$ such that $\gamma(q_1)$ and
$\gamma(q_2)$ do not lie on
$\Gamma^\bullet_{q_1,q_2}$.
Then as in the first paragraph of the proof of Lemma 4.2,
one sees that the conic $\Gamma^\bullet_{q_1,q_2}$
does not pass through $\gamma(a)$. Hence we may assume that
the conic conic $\Gamma^\bullet_{p_1,p_2}$ itself does not
pass through $\gamma(a)$.

Let $I$ be a closed interval in $(0,2\pi)$ such that
 $p_1,p_2\in I$.
Then the function $f_2^\bullet$ defined in Example 2.6 (ii)
is an intrinsic system system of order $5$ on $I$ satisfying
the boundary regularity condition $(5,5)$.
We can prove the existence of a sextactic  point
on $(0,2\pi)$ with the property that
the osculating conic at $s$ is inscribed
and does not meet $\gamma(a)$ with methods as in the proof of Lemma 4.2;
see also Lemma 4.10 in [TU2] where this is also proved.

If furthermore $\gamma$ is $C^1$ at $a$, one can also assume that a
minimal circumscribed conic $\Gamma^\circ_{p_1,p_2}$ does not pass
through $\gamma(a)$ and show the existence of sextactic  point
on $(0,2\pi)$ with the property that
the osculating conic at $s$ is circumscribed
and does not meet $\gamma(a)$. Here we use the fact that
the function $f_2^\circ$ defined in  Example 2.6 (ii)
is an intrinsic system system of order $5$ on $[\varepsilon,
2\pi-\varepsilon]$ satisfying the boundary regularity condition $(5,5)$.

It is now straightforward how the arguments in
the proof of the Jackson Lemma for Flexes of Functions 4.1 carry over
to the present situation. We only sketch the main points. Let
us assume that we are dealing
with circumscribed conics assuming that $\gamma$ is $C^1$ at $a$.
(The existence of a clean sextactic point whose osculating conic is
inscribed
is similar except that we do not need to assume $\gamma$
to be $C^1$ at $a$ since we have already proved the existence of sextactic
point  
on $(0,2\pi)$ with the property that the osculating conic at $s$
is circumscribed and does not meet $\gamma(a)$.)

Let $\Phi(\gamma)$ be the set of sextactic points such that the osculating
function there are circumscribed and do not meet $a$. We have already seen
that $\Phi(\gamma)$
is nonempty. For $p\in \Phi(\gamma)$ we define $I(p)$ to be the smallest
closed interval containing
$F_2^\circ(p^2)=\Gamma^\circ_{p^2}\cap\gamma$, but not containing $a$.
Analogous
to $u_p$ one defines $\gamma_p$
as $\gamma$ on the interval $I(p)$ and equal to $\Gamma_{p^2}^\circ$ on the
complement of $I(p)$. Notice
that $\gamma_p$ is a closed simple contractible curve in $P_2$. The partial
ordering on
$\Phi(\gamma)$ is now defined by setting
$p<\!\!<q$ for 
$p,q\in
\Phi(u)$ if
$$
I(p)\subset I(q) \qquad {\rm and } \qquad  D(\gamma_q)\subset D(\gamma_p),
$$
where $D(c)$ denotes the closed contractible region bounded by a simple
closed contractible curve $c$. The
same arguments as in the proof of Theorem 4.1 can now be used to show that a
point $p\in\Phi(\gamma)$ that is
minimal with respect to this ordering is a clean sextactic point of the type
we are are trying to find. The
final step is again to use Zorn's Lemma to prove the existence of a minimal
point in $\Phi(\gamma)$. Let $S$
be a totally ordered subset of $\Phi(\gamma)$. Set
$$
D={\bigcup_{p\in S}}D(\gamma_p).
$$
The boundary of $D$ consists of two pieces: One piece is the image of the
interval 
$$I=\bigcap_{p\in S}I(p)$$
under $\gamma$. We would like to show that the other piece, the complement
of $\gamma(I)$ which we denote by
$\Gamma$, is an arc of a conic. We choose five different points
$q_1,\dots,q_5$ on $\Gamma$ and five
sequences $(q_{1,k}),\dots, (q_{5,k})$ of points on a sequence  conics
$\Gamma^\circ_{p_k^2}$ such that
$(q_{i,k})$ converges to $q_i$ for $i=1,\dots,5$. We can assume that the
corresponding sequence $(p_k)$ converges to a
point $p_\infty$. We can now use arguments as in the proof of Theorem 4.1 to
show that arcs of
$\Gamma^\circ_{p_k^2}$ converge to $\Gamma$ and that $\Gamma$ is an arc of
the osculating conic
$\Gamma^\circ_{p_\infty^2}$. It follows that $p_\infty\in\Phi(\gamma)$ and
$p_\infty<\!\!<p$ for all $p\in
S$ and we can apply Zorn's Lemma.
\qed

\bigskip
\noindent
{\bf \S 5 \hskip 0.1in  On the existence of $2n+2$ clean flexes on
periodic functions}
\medskip
We proved in Section three that a smooth periodic function $u$ has at
least three flexes where the osculating functions are greater or equal to
$u$
and similarly at least three flexes where the osculating functions are less
than or equal
to $u$. In Section four we also started to study the existence of clean
flexes in the
Jackson Lemma and this
will be continued in this section.

Our main result here is the following theorem which is the same
as Theorem 1.1 in the introduction if ${\cal A}=A_{2n+1}$.
 \medskip 
 
\ni{\bf
5.1. Theorem.} {\it  Let ${\cal A}$ be a Chebyshev space of order
$2n+1$ where $n\ge 1$.
Let $u$ be a $C^{2n}$-function on $S^1$ which does not
belongs to ${\cal A}$.
Then $u(t)$ has at least  $n+1$ different (intervals of) clean maximal
${\cal A}$-flexes
and at least $n+1$ different (intervals of) clean minimal ${\cal
A}$-flexes.}
\medskip

\ni{\bf Remark.} Notice that a clean maximal ${\cal A}$-flex
cannot be a clean minimal ${\cal A}$-flex if $u$ does not
belong to ${\cal A}$.
The theorem therefore gives us $2n+2$ clean flexes.
\medskip

\ni{\bf Example.} Theorem 5.1 is optimal.
Set ${\cal A}=A_{2n+1}$. The flexes of a given
function $u$ are the zeros of $L_{2n+1}u$ where $L_{2n+1}$ is the operator
$$
L_{2n+1}=D(D^2-1)\cdots (D^2-n^2),
$$
where $D={d /dt}$,
see Proposition A.7 in Appendix A.
If we set $u(t)=\sin(n+1)t$, then $L_{2n+1}u(t)$ is
proportional to $u(t)$. Thus $u(t)$
has exactly $2n+2$-flexes which are all clean.
\medskip

The main new tool in the proof of Theorem 5.1 is the
following proposition, which is in
the same spirit as the induction argument in [Ba], p.~201--204; see also
Hilfssatz 2 on page 140 in [N\"o].
Notice that unlike here, global properties
are not treated in [Ba] and [N\"o].

\medskip
\ni{\bf 5.2. Proposition.}
{\it 
Let ${\cal A}$ be a Chebyshev space of order
$2n+1$ where $n\ge 2$.
Let $u$ be a $C^{2n}$-function on $S^1$ and $f_u$ the
corresponding intrinsic system of order $2n+1$ on $S^1$.
\item{(i)} Suppose that
$$ 
f_u(p,a^{n-1})(a)\ge 2n,
$$
and that $F_u(p,a^{n-1})\setminus \{p\}$
is a closed interval
where $p$ and $a$ are two different points.
Then both $(p,a)$ and $(a,p)$ contain a clean maximal ${\cal A}$-flex.
\item{(ii)} Let
$[a,b]$ be a nontrivial closed interval on $S^1$.
Suppose that
$$ f_u(p,a^{n-1})(a)\ge 2n \qquad{\rm and}\qquad
 f_u(p,b^{n-1})(b)\ge 2n
$$
for a point $p\not\in [a,b]$.
If furthermore both
 $F_u(p,a^{n-1})\setminus \{p\}$
and $F_u(p,b^{n-1})\setminus \{p\}$
are different closed intervals, then there is
a clean
maximal ${\cal A}$-flex in the interval $(a,b)$.
\item{}}

\medskip
\ni{\it Proof.}  We first proof the claim in (i).
 Instead of the original function $u$ we consider the following
function (See Figure 5.)
$$
v(t)=\cases{
 u(t) &
\quad {\rm if $t\in [p,a]$}, \cr
\varphi_{u,(p,a^{n-1})}(t)=\varphi_{u,(a^{n})}(t) & \quad otherwise, \cr
} 
$$
where $\varphi_{u,(p,a^{n-1})}$ is a minimal function
with respect to $(p,a^{n-1})$.

\medskip
\centerline{\epsfxsize=3in \epsfbox{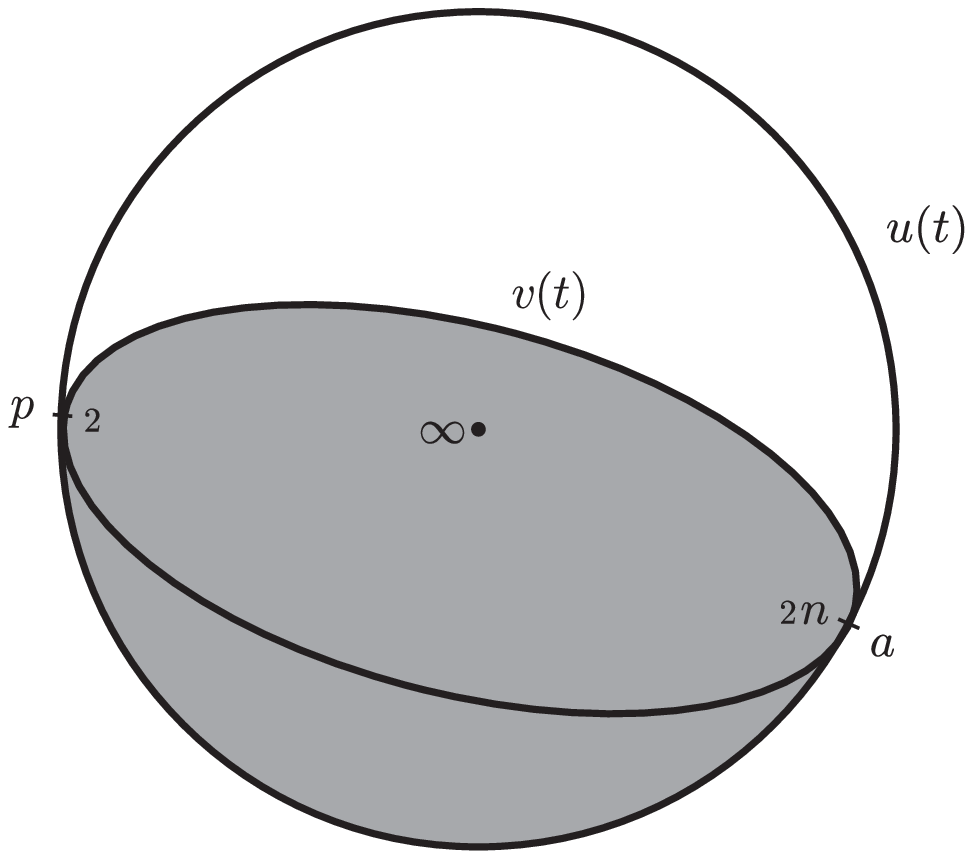}}
\centerline{Figure 5.}
\bigskip

\noindent
The function $v(t)$ is at least $C^{1}$ in $p$
and  at least $C^{2n-1}$ in $a$. Set $I=[p,a]$.
Then $v$ satisfies the boundary
condition $(1,n)$ on $I$; see Definition 2.2.
We consider the intrinsic system $f_v$ defined on $I^n_{(1,n)}$.
We have that 
$$
f_v(p,a^{n-1})(p)+f_v(p,a^{n-1})(a)\ge 2(n+1)
$$
We can therefore apply the Contraction Lemma 3.9 to this situation
which now implies that there is an $f_v$-flex $s$ in the
open interval $(p,a)$ which is clearly also an $f_u$-flex.
We have nothing to prove if $s$ is a clean flex of $u$. Therefore
we assume that it is not. The osculating conic $\varphi_{u,(s^n)}$
can clearly only take on the same values as $u$  in the interval
$[p,a]\cup F_u(p,a^{n-1})$ since $v$ is strictly larger than $u$
on its complement. Let $q\in F_u(s^n)$ be point that is not in the
same component of $F_u(s^n)$ as $s$. Let us assume that $q$ comes
before $s$ in the interval $[p,a]\cup F_u(p,a^{n-1})$, the other case being
similar. Define a function $w$ by setting
$$
w(t)=\cases{
u(t)  &\quad {\rm if $t\in [q,s]$}, \cr
\varphi_{u,(s^{n})}(t) & \quad otherwise. \cr
} 
$$
This function $w$ is $C^{2n}$ in $s$ since $s$ is a flex. Hence $w$ is
$C^{2n}$ except possibly in $q$ where it is at least $C^1$. We
can
now apply the Jackson Lemma 4.1 to $w$. It follows that $w$ has a clean flex
$s'$ whose osculating function does not take on the same value as $u$ in
$q$. Hence $s'$ must be contained in the interval $(q,s)$. Notice also that
$s'$ 
cannot be contained in the interval $F_u(p,a^{n-1})\setminus \{p\}$ since
then  $\varphi_{u,(s')}=\varphi_{u,(p,a^{n-1})}$ contradicting that $s'$ is
a clean 
flex. It follows from this discussion that $s'\in(p,a)$.

The proof that $(a,p)$ contains a clean flex is very similar.

Next we prove (ii).  Instead of the original function $u(t)$, we consider
the following 
function (See Figre 6.)
$$
v(t)=\cases{
\varphi_{u,(p,a^{n-1})}(t)=\varphi_{u,(a^{n})}(t)  &
\quad {\rm if $t\in [p,a]$}, \cr
u(t) & \quad if $t\in [a,b]$, \cr
\varphi_{u,(p,b^{n-1})}(t)=\varphi_{u,(a^{n})}(t)  & \quad
{\rm if $t\in
[b,p+2\pi]$}, \cr} 
$$
where $\varphi_{u,(p,a^{n-1})}$ (resp.~$\varphi_{u,(p,b^{n-1})}$)
are minimal function
with respect to $(p,a^{n-1})$ (resp. $(p,b^{n-1})$).

\medskip
\centerline{\epsfxsize=3in \epsfbox{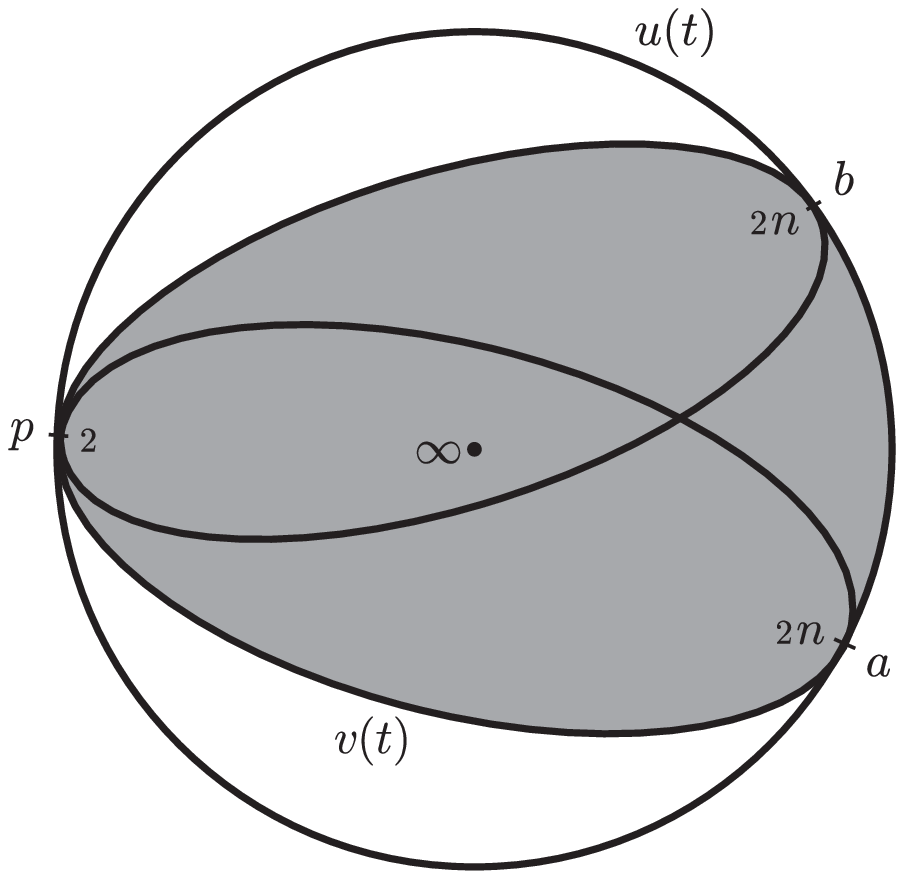}}
\centerline{Figure 6.}
\bigskip

\noindent

The function $v(t)$ is at least $C^{1}$ in $p$
and  $C^{2n-1}$ in $a$ and $b$.
Moreover $v$ 
 satisfies at least the boundary regularity condition $(n,n)$ on $[a,b]$;
see Definition 2.2.
Let $f_v$ denote the intrinsic system
of order $2n+1$ on  $[a,b]$
satisfying the boundary regularity condition $(n,n)$
whose existence was proved in Proposition 2.4.
We choose  $n$ different points
$
p_1\prec \cdots \prec p_{n}
$
in the interval $(a,b)$, but not in $F_u(p,a^{n-1})\cup
F_u(p,b^{n-1})$.
Arguing exactly as at the beginning of the proof of the
Jackson Lemma 4.1, we can assume that
$$
f_v(p_1,\cdots , p_{n})(p)= 0.
$$
We choose a point $q\in (p,p_n)$ as follows:
First we consider the case
$$
\sum_{i=1}^n
 f_v(p_1,\cdots , p_{n})(p_i)=2n.
$$
By Axiom (A6), there exists a
point $q \in F_v(p_1,\cdots , p_{n})$
that is different from $p$ and $p_1,\dots,p_{n}$.
Without loss of generality, we may assume that $q\in (p,b)$.
After interchanging $q$ and $p_1,\dots,p_n$ if necessary,
we may assume 
$$
q\prec p_1\prec \dots \prec p_n.
$$
Next we consider the case that
$$
\sum_{i=1}^{n}
f_v(p_1,\cdots , p_{n})(p_i)\ge 2n+2.
$$
In this case we set $q=p_1$.
In both of these two cases,
there is a sufficiently small $\varepsilon>0$,
such that $p+\varepsilon\prec q$ and
$p_{n-1} \prec p_{n}-\varepsilon$.
Moreover, $v$ is $C^{2n-1}$ on $[p+\varepsilon,p_{n}-\varepsilon]$
satisfying at least the boundary regular condition $(n-1,n-1)$.
By Proposition 2.8, we can associate to $v$ an
intrinsic system $f'_v$ of
order $2n-1$ on the interval $[p+\varepsilon,p_{n}-\varepsilon]$
with base point $p_{n}$.
We have that
$$
\sum_{t\in [q,p_{n-1}]} f'_v(p_1,\dots,p_{n-1})(t)\ge 2n.
$$
Hence there exists a point $r\in (q,p_{n-1})$ by the
Contraction Lemma 3.10 such that
$$
 f'_v(r^{n-1})(r)\ge 2n.
$$
If $r\not\in [a,b)$, then $r\in (p,a)$.  In this case
$\varphi_{v, (p_{n},r^{n-1})}$ is locally around $r$ greater or equal to
$\varphi_{u,(a^{n})}$
and the two functions meet in $r$
with multiplicity at least $2n$.
Here
$\varphi_{v, (p_{n},r^{n-1})}$
is the minimal function of $v$ with respect to
$(p_{n},r^{n-1})$.
On the other hand, the value of $\varphi_{v, (p_{n},r^{n-1})}$ in the point
 $p_n$ is smaller than that of
 $\varphi_{u,(a^{n})}$ in $p_n$.
This means that $\varphi_{v,(p_{n},r^{n-1})}$ meets
$\varphi_{u,(a^{n})}$
with multiplicity at least $2n+2$
which implies that
$\varphi_{v,(p_{n},r^{n-1})}=\varphi_{u,(a^{n})}$.
Hence
$p_{n}\in F_u(a^{n})$ which is a contradiction.
Thus we have $r\in [a,b)$.

We have
$$
 f_v(r^{n-1},p_n)(r)+ f_v(r^{n-1},p_n)(p_n)\ge 2(n+1).
$$
Hence we can use the Contraction Lemma 3.9 (or 3.10) to prove
the existence of a
$f_v$-flex $s$ in the open interval $(a,p_n)\subset(a,b)$.
The point $s$ is clearly also an $f_u$-flex, but might not
be a clean flex. We can now as in the proof of part (i) of this proposition
introduce a function $w$ and apply the Jackson Lemma to show
the existence of a clean flex that must be contained in the interval
$(a,b)$.  
\qed

\medskip
\noindent
{\it Proof of Theorem 5.1.}
We shall now prove the existence of $n+1$ clean maximal ${\cal A}$-flexes on
$u$
by induction over the order $2n+1$ of the Chebyshev space ${\cal A}$.
The result can then be applied to $-u$ to also prove the existence
of $n+1$ clean maximal ${\cal A}$-flexes on $u$.

If $n=1$, it follows quite easily from the proof of the Kneser
Lemma 3.6 that $u$ has two clean maximal ${\cal A}$-flexes.
We now assume the claim of the theorem to be true for $n-1$ and show that
it then follows for $n$.
We fix $p\in S^1$ that is not a flex.
By Theorem A.2 in Appendix A there exists a function
$\psi_0\in {\cal A}$ such that
$$ 
u(p)=\psi_0(p),\qquad u'(p)=\psi'_0(p).
$$
Now we define a linear subspace $V$ of ${\cal A}$ by setting
$$
V=\left\{ \varphi \in {\cal A}\,;\,
\varphi(p)=0,\,\, \varphi'(p)=0\right\}.
$$
We set
$$
\psi_1(t)=\alpha + \beta \cos t + \gamma \sin t
$$
We can adjust the three coefficients $\alpha,\beta,\gamma$
so that $\psi_1$ satisfies
$$
\psi_1(p)=\psi'_1(p)=0.
$$
Since $\psi_1$ belongs to $A_{3}$, it has at most
two zeros counted with multiplicities.
Thus $\psi_1$ has no zeros other than $p$
and the second derivative of $\psi_1$ at $p$
does not vanish. We then set
$$
{\cal A}_{\psi_1}
=\left\{ {\varphi \over \psi_1}
\,;\,\varphi\in V
\right\}.
$$
It can easily be checked that ${\cal A}_{\psi_1}$ is a Chebyshev
space of order $2n-1$ since ${\varphi/\psi_1}$
is $C^{2n-2}$ at $p$; see Appendix B. We set
$$
v(t)={u-\psi_0\over \psi_1}.
$$
Then $v$ is a $C^{2n-2}$-function on $S^1$.
So by applying the induction assumption
there exist $n$ clean ${\cal A}_{\psi_1}$-flexes
$s_1,\dots,s_n$.
Let $\varphi_j$ in ${\cal A}_{\psi_1}$  be the
osculating function of $v$ at $s_j$.
Since $\varphi_j$ is a function in  ${\cal A}_{\psi_1}$,
there exists a function $\hat \varphi_j\in V$ such
that $\varphi_j=\hat \varphi_j/\psi_1$.
If some $s_j$ is equal to $p$, then
$\hat \varphi_j$ meets $u$ only in one component,
and this implies that $p$ is a clean ${\cal A}$-flex, which is a
contradiction.

So none of the points $s_j$ can be equal to $p$.
Hence $\psi_1$ does not vanish in any of the points $s_j$
and it follows that the first $2n-1$ derivatives
of $(u-\psi_0)-\hat\varphi_j$  vanish in $s_j$.
Since $(u-\psi_0)-\hat\varphi_j$ vanishes with multiplicity at
least two in $p$, we have that
$$
\psi_0+\hat \varphi_j=\varphi_{(p,s_j^{n-1})}=\varphi_{(s_j^n)},
$$
where $\varphi_{(p,s_j^{n-1})}$ and $\varphi_{(s_j^{n})}$
are the maximal functions of $u$ with respect to the
$n$-uples $(p,(s_j)^{n-1})$ and $((s_j)^{n})$ respectively.
This implies that
$$
f_u(p,s_j^{n-1})(s_j)\ge 2n, \quad f_u(p,s_{j-1}^{n-1})(s_{j-1})\ge 2n
\qquad {\rm for}\qquad j=0,1,\dots,n,n+1,
$$
where we have set $s_0=s_{n+1}=p$.
Since each $s_j$ is a clean ${\cal A}_{\psi_1}$-flex,
$F_u(p,s_j^{n-1})\setminus \{p\}$ is a closed interval for
$j=1,...,n$.
By Proposition 5.2, there is a clean maximal ${\cal A}$-flex $t_j$
on each of the intervals $(s_{j-1},s_j)$ for $j=0,\dots,n+1$.
\qed

\medskip

We now also prove the other theorem stated in the introduction.
The proof is similar to the one of Theorem 5.1.
 
\medskip
\noindent
{\bf 5.3 Theorem.}
{\it
Let $\gamma:S^1\to {\bf R}^2$ be a strictly convex
$C^4$-regular curve which is not a conic.
Then $\gamma$ has at least three (intervals of) clean
maximal and at least three (intervals of) clean
minimal  sextactic points.}\medskip

\noindent
{\it Proof.} We only prove the existence of three clean maximal
sextactic points, the proof of the existence of three clean
 minimal flexes being similar.
 Let $f_2^\bullet$ be the
 intrinsic system of order $5$ on $S^1$ introduced in Example 2.6 (ii).
The following lemma is analogous to Proposition 5.2 (i) and has a very
similar
proof, which we therefore omit only remarking that it
is this time based on the Jackson Lemma for Sextactic
Points 4.3.

\medskip
\ni{\bf 5.4. Lemma.}
{\it 
Suppose that
$$ 
f_2^\bullet(p,a)(a)\ge 4,
$$
and that $F^\bullet(p,a)\setminus \{p\}$
is a closed interval where $p$ and $a$ are two different points.
Then both $(p,a)$ and $(a,p)$ contain a clean maximal sextactic point. \qed
}

\medskip
We now come back to the proof of Theorem 5.3.
By Theorem 4.3 there is a clean maximal sextactic point $p$ on the curve
$\gamma$.
 We can also use
Theorem 4.3 to
find a clean maximal sextactic point $q$ whose osculating conic does not
meet 
$p$.
We first show that $f_2^\bullet(p,q)(p)=2$ and $f_2^\bullet(p,q)(q)=2$.
Since the osculating conics coincide with the maximal inscribed
conic at $p$ and $q$, the inequality $f_2^\bullet(p,q)(p)\ge 4$
(or $f_2^\bullet(p,q)(q)\ge 4$) implies
 $f_2^\bullet(p,q)(p)=f_2^\bullet(p^2)(p)=\infty$ (or
 $f_2^\bullet(p,q)(q)=f_2^\bullet(q^2)(q)=\infty$) and
 the osculating conic at
$p$ would pass through $q$ (or the one at $q$ through $p$.) This is a
contradiction.
Hence $f_2^\bullet(p,q)(p)=f_2^\bullet(p,q)(q)=2$.

By (A6), there is a point $r$ distinct from $p,q$
such that $f^\bullet(p,q)(r)\ge 2$. Assume that $p\prec q\prec r$.
The restriction of $f^\bullet_2$ to $[p,r]$
is an intrinsic system of order $5$
satisfying the boundary regularity condition $(2,2)$.
We can therefore define $f_{2,p}^\bullet$ as before Lemma 3.7., i.e.,
we set for $x\in [p,r]$ and $y\in S^1$
$$
f_{2,p}^\bullet(x)(y)=\cases{f^\bullet_2(p,x)(y) &  if $y\not =p$,\cr
                 f^\bullet_2(p,x)(y)-2 &  if  $y=p$. \cr}
$$
By Lemma 3.7 we know that $f_{2,p}^\bullet$ restricted to $[q,r]$ is an
intrinsic 
system of order $3$  satisfying at least the boundary regularity
condition $(1,1)$. There is now a subinterval $[a,b]$ of $[q,r]$
satisfying the conditions in the Kneser Lemma 3.6 which implies that there
is
a point $c$ in $(a,b)$ such that $F_{2,p}^\bullet(c)$ is connected. This
implies 
that $f_2^\bullet(p,c)(c)\ge 4$
and that $F_2^\bullet(p,c)\setminus \{p\}$ is a closed interval. By Lemma
5.4 we then have 
a clean maximal point in the interval $(p,c)$ and another one in $(c,p)$.
(One of these
sextactic points might coincide with $q$.)
Since $p$ is also a clean sextactic point we have proved the existence of
thee (intervals of)
clean maximal sextactic points.
 \qed

\bigskip
\noindent
{\bf \S 6 \hskip 0.1in  Arrangements of clean flexes.}
\medskip
As was pointed out in Example 2.6 (i),
there are two clean maximal vertices and
two clean minimal vertices on a given simple closed
curve $\gamma:S^1\to {\bf R}^2$. It is a natural question to ask
in which order the clean maximal and the clean minimal vertices are arranged
on $S^1$. In [TU2], the authors proved that
there are four points $t_1\prec t_2\prec
t_3\prec t_4$ on $S^1$
such that $t_1,t_3$ are clean maximal vertices and
$t_2,t_4$ are clean minimal vertices.
Now that we have proved the
existence of $2n+2$ clean ${\cal A}$-flexes on
a $2\pi$-periodic function $u$, we can ask again how the maximal and minimal
ones
are arranged relative to each other.
We will say that the {\it clean ${\cal A}$-flexes on $u$ change sign at
least
$m$-times}
if there are $2m$-points
$$
p_1\prec q_1\prec \cdots\prec p_m \prec q_m
$$
on $[0,2\pi)$ such that $p_j$ for $j=1,\dots,m$ are
clean maximal ${\cal A}$-flexes and $q_j$ for $j=1,\dots,m$ are
clean minimal ${\cal A}$-flexes.

\medskip
\ni{\bf 6.1 Theorem.} {\it Let $u$ be a $2\pi$-periodic
$C^{2n}$-function which is not in $\cal A$. Then the clean
${\cal A}$-flexes on $u$
change sign at least four times.}

\medskip
We do not know whether Theorem 6.1 gives an optimal lower bound
on the number of sign changes or not.

To prove the theorem, we will use the
abstract theory of pairs of intrinsic circle systems of which we give a
quick review.

\medskip
\noindent
{\bf 6.2 Definition.}
A family of nonempty closed subsets $F=(F_p)_{p\in S^1}$ of
$S^1$
is called an {\it intrinsic circle system} on $S^1$
if it satisfies the following three conditions for any $p \in S^1$.
\item{(I1)}
 If $q\in F_p$, then $F_p=F_q$.
\item{(I2)} 
If $p'\in F_p$, $q'\in F_q$ and
$q\succeq p'\succeq q' \succeq p (\succeq q$),
then $F_p=F_q$ holds.
\item{(I3)} Let $(p_n)_{n\in {\bf N}}$
and $(q_n)_{n\in {\bf N}}$ be two sequences in $S^1$ such that
      $\lim_{n\to \infty}p_n=p$ and $\lim_{n\to \infty}q_n=q$
      respectively.
      Suppose that $q_n \in F_{p_n}$ for all $n\in {\bf N}$. Then
$q\in F_p$ holds.

\medskip
A pair of intrinsic circle systems
$(F^+, F^-)$ is said to be {\it compatible}
if it satisfies the following two conditions.
\item{(C1)}
$F^+_p\cap F^-_p=\{p\}$
for all $p\in S^1$.
\item{(C2)} Suppose that $F^+(p)$
(resp.~$F^-(p)$) is connected.
Then there are no points $q$
in a sufficiently small
neighborhood of $p$
such that 
$F^+(q)$ 
(resp.~$F^-(q)$) is connected.

\medskip
In [TU2], Theorem 1.4, the authors proved the following

\medskip
\ni{\bf 6.3 Lemma.} {\it
Let $(F^+, F^-)$ be a compatible pair of
intrinsic circle systems.
Then 
there are four points $p_1,p_2,p_3,p_4\in S^1$ satisfying
$p_1\succ p_2 \succ p_3 \succ p_4 (\succ p_1)$
such that $F^+_{p_1}, F^+_{p_3}$
and $F^-_{p_2}, F^-_{p_4}$ are connected subsets of $S^1$.
}

\medskip
We now prove the theorem as a corollary of Lemma 6.3.
\medskip

\noindent
{\it Proof of Theorem 6.1.}
Suppose that the claim in the theorem is not true.
Then the clean flexes on $u$
change sign exactly two times.
There are clearly disjoint closed intervals $I$ and $J$
containing all the clean maximal
and all the clean minimal flexes respectively.  There is a
point 
 $p\in S^1$ such that $p\not\in I\cup J$.
Without loss of generality, we may set
$p=0$, and 
$$
0<{\rm inf}(I)<{\rm sup}(I)<
{\rm inf}(J)<{\rm sup}(J)<2\pi.
$$
We set 
$$
\eqalign{
F^+_q&=\cases{F_{u}(q,p^{n-1})\setminus \{p\} &  if  $q\not \in F(p^n)$
\cr F_{u}(p^{n}) &  if  $q\in F(p^n)$, \cr } \cr
F^-_q&=\cases{F_{-u}(q,p^{n-1})\setminus \{p\} &  if  $q\not \in F(p^n)$
\cr F_{-u}(p^{n}) &  if  $q\in F(p^n)$. \cr }
}
$$
It is easy to check that $(F^+,F^-)$ is
a compatible pair of intrinsic circle systems.
By Lemma 6.3, there are four points $p_1,p_2,p_3,p_4\in S^1$ satisfying
$p_1\succ p_2 \succ p_3 \succ p_4 (\succ p_1)$ such that $F^+_{p_1},
F^+_{p_3}$ and $F^-_{p_2}, F^-_{p_4}$ are all connected.
Now we claim that $p_1,p_3\in I$.
Assume that $p_1\not\in I$ holds, then
we have
$p_1 \in (0,\inf(I))$ or $p_1 \in (\sup(I),2\pi)$.
Since  $F^+_{p_1}$ is connected, we have
$$
f_u(p_1,p^{n-1})(p_1)\ge 4,\qquad
f_u(p_1,p^{n-1})(p)\ge 2n-2.
$$
By Lemma 3.7, we can define a new intrinsic system $g$ of
order $2n-2$ on $[0,p_1]$ (resp. on $[p_1,2\pi]$) by
$$
g(q_1,...,q_{n-1})(r)=\cases{f(p_1,q_1,\dots,q_{n-1})(r) &  if  $q\not
=
r$,\cr
                 f(p_1,q_1,\dots,q_{n-1})(r)-2 &  if  $r=q$, \cr},
$$
where $(q_1,...,q_{n-1})\in [0,p_1]^{n-1}_{(1,n)}$
(resp. $(q_1,...,q_{n-1})\in [p_1,2\pi]^{n-1}_{(1,n)}$).
Since $g$ satisfies the boundary regularity condition $(1,n)$
on $[0,p_1]$ and $(n,1)$ on $[p_1,2\pi]$, there
are  clean $g$-flexes  $s_1\in [0,p_1]$ and
$s_2\in [p_1,2\pi]$. (Apply Proposition 5.2 (i) to 
 the Chebyshev space $A_{\psi_1}$ defined in the
proof of Theorem 5.3.) Thus we have
$$
f(p_1,s_j^{n-1})(p_1)\ge 2,\qquad f(p_1,s_j^{n-1})(s_j)\ge 2n \qquad
(j=1,2).
$$
By  Proposition 5.2 (i),
 we have a clean maximal $\cal A$-flex $s'_1\in (0,s_1)$ and $s'_2\in
 (s_2,2\pi)$ respectively.
If $p_1\in (0,\inf(I))$ (resp. $p_1\in (\sup(I),2\pi)$), we have
$s'_1\in (0,\inf(I))$ (resp. $s'_2\in (\sup(I),2\pi)$).
This is a contradiction since all of the clean maximal flexes
are contained in $I$. So we have $p_1\in I$.
Similarly we can see that $p_3\in I$
and $p_2,p_4\in J$.
This contradicts the relation $p_1\succ p_2\succ p_3\succ p_4$.
Hence we must have at least four sign changes of clean flexes.
\qed

Finally we formulate two open problems
which are in our opinion the most important ones on flexes of periodic
functions.

\medskip
\noindent
{\bf Problems.} (1) Give a best lower bounds for
the number of sign changes of clean ${\cal A}$-flexes of a periodic
function.
The number $2n+2$ is a tempting guess.

(2)  Let $u$ be a  $2\pi$-periodic
$C^\infty$-function and ${\cal S}$ be the union over the sets of
$A_{2n+1}$-flexes
on $u$ where $n$ ranges over all natural numbers. Is the set $\cal S$ a
dense subset of $S^1$?

\bigskip
\noindent
{\bf Appendix A. \hskip 0.1in  Chebyshev spaces.}
\medskip

In this appendix, we shall define Chebyshev spaces as certain
linear subspaces of smooth functions, bring their basic theory and
explain the existence of $2n+2$ flexes on a periodic
function as an application.

Chebyshev spaces are related to Fourier series and disconjugate
operators.
It is well known that many theorems in the spirit of
the classical four vertex theorem can be proved
using either Fourier series (see Hayashi [Hy]
and Blaschke [Bl], p.~68) or disconjugate operators
(see [Ar2], [GMO], [OT] and [Ta]).

A $C^n$-function $u$ defined on  ${\bf R}$ has a {\it zero of
order $m$} (or a {\it zero with multiplicity $m$})
at $s$ where $1\le m\le n$
if the value and the first $m-1$ derivatives of $u$ at $s$ vanish,
but not the $m$-th. Notice that we do not define the order of a zero $s$ of
a $C^n$-function $u$ that vanishes in $s$ together with all its $n$
derivatives since it will
not be needed in the following.

\medskip
\noindent
{\bf A.1 Definition.} A linear subspace
${\cal A}$ of $C^{n-1}({\bf R})$ is called
a {\it  Chebyshev space of order $n$ } if the following conditions are
satisfied.
\item{(i)} $n\le {\rm dim}\, {\cal A}$.
\item {(ii)} Every nonvanishing function $\varphi$ in ${\cal A}$
 has at most $n-1$ zeros counted with multiplicity in ${\bf R}/2\pi {\bf
Z}$.
\item{(iii)} Every function $\varphi$ in ${\cal A}$ is
 $2\pi$-periodic if $n$ is odd.
\item{(ii)}Every function $\varphi$ in ${\cal A}$
is  $2\pi$-antiperiodic, i.e., $\varphi(t+2\pi)
=-\varphi(t)$, if $n$ is even.

\medskip
In the definition of a Chebyshev space we allow the possibility
that ${\rm dim}{\cal A}=\infty$. We will see in Theorem A.2 that
the dimension of a Chebyshev space is always
finite and equal to its order $n$.

 We now give a few examples of Chebyshev spaces.

The space 
$$
A_{2k+1}=
\left\{\varphi\in C^{2k}({\bf R})\,;\,
\varphi(t)=a_0 +\sum_{n=1}^k (a_n \cos nt + b_n \sin nt) \right\}.
$$
is a Chebyshev space of order $2k+1$.
In fact, we have
$$
\cos nt={z^n+z^{-n}\over 2},\qquad \sin nt={z^n-z^{-n} \over 2i},
$$
where $z=e^{it}$.
The functions $z^k \varphi$ for $\varphi$ in $A_{2k+1}$ are
polynomials in $z$ of degree less than or equal to $2k$.
Consequently, the number of zeros of the functions $\varphi$ in $A_{2k+1}$
can at most be
$2k$. 

Similarly, the space
$$
A_{2k}=
\left\{\varphi\in C^{2k}({\bf R})\,;\,
\varphi(t)=\sum_{n=1}^k \left[
a_n \cos \left({(2n-1)t \over 2}\right)
 + b_n \sin \left({(2n-1)t\over 2}\right)\right] \right\}.
$$
is a Chebyshev space of order $2k$.

A linear differential operator $L$ of order $n$ on ${\bf R}$
is called {\it disconjugate} if its kernel ${\rm Ker} L$
is a Chebyshev space of order $n$.
Examples of disconjugate operators are
$$
\eqalign{
L_{2k+1}&=D(D^2-1)\cdots (D^2-k^2), \cr
L_{2k}&=(D^2-\left({1\over 2}\right)^2)
\cdots (D^2-\left({2k-1\over 2}\right)^2),
}
$$
where $D={d /dt}$, 
since their kernels are $A_{2k+1}$ and $A_{2k}$ respectively. We refer to
Proposition A.6 for more
information on Chebyshev spaces and disconjugate operators.

Further simple examples of Chebyshev spaces can be obtained as follows. Let
$\cal A$
be the linear span of either one of the sets
$$\eqalign{
A_{2k+1}^n&=
\left\{\varphi_1\cdots \varphi_n
\,;\, \varphi_1,\dots, \varphi_n\in A_{2k+1} \right\}, \cr
A_{2k}^n&=
\left\{\psi_1\cdots \psi_n
\,;\, \psi_1,\dots, \psi_n\in A_{2k} \right\}.
}
$$
Then ${\cal A}$ is a Chebyshev space of order
$n(2k+1)$ in the first case and of order $2nk$ in the second case.

The following property of Chebyshev spaces is crucial.

\medskip
\noindent
{\bf A.2 Theorem.} {\it Let ${\cal A}$ be a
Chebyshev space in 
$C^{n-1}({\bf R})$ of order $n$.  Let
$$
0\le t_1\le  t_2\le  \cdots \le  t_{n} <2\pi
$$
be $n$ points and $\nu_j$ the multiplicity with
which $t_j$ occurs as a component of
the $(j-1)$-uple $(t_1,\dots,t_{j-1})$.
Then the linear map $T:{\cal A}\to {\bf R}^n$ defined by
$$
 T(\varphi)=(\varphi^{(\nu_1)}(t_1),\dots,\varphi^{(\nu_n)}(t_n))
$$
is an isomorphism. In particular,  ${\cal A}$ is finite dimensional
and its dimension coincides with
the order $n$.}

\medskip
\noindent
{\it Proof.} Let $\varphi\in {\cal A}$ be in the kernel of the map $T$. It
follows from 
the definition of $T$ that $\varphi$ vanishes at least $n$ times counted
with
multiplicities. The definition of a Chebyshev space now implies that
$\varphi$ vanishes identically.
We have therefore proved that $T$ is injective and it follows that
$\dim\;{\cal A}\le n$.
By the definition of a Chebyshev space we have $\dim\;{\cal A}\ge n$. Hence
$\dim\;{\cal A}=n$ and $T$ is
an isomorphism. \qed

\medskip
\noindent
{\bf A.3 Corollary.} {\it Let ${\cal A}$ be a
Chebyshev space  of order $n$.  Let
$$
0\le t_1< t_2< \cdots < t_{k} <2\pi
$$
be $k$ different points where $k\le n-1$ and let
$
i_1, i_2,\dots,i_k
$
be $k$ positive integers satisfying
$i_1 + i_2 + \cdots + i_k=n-\ell$ for $\ell\ge0$.
Then the set of functions $\varphi$ in ${\cal A}$
which have  zeros of order $i_j$ at $t_j$ for $j=1,\dots,k$ is an
$\ell$-dimensional 
subspace of $\cal A$.}\qed

We shall now prove the following

\medskip

\noindent
{\bf A.4 Theorem.} {\it Let ${\cal A}$ be a
Chebyshev space in $C^{n-1}({\bf R})$ of odd
(resp. even) order $n$.
Let $u$ be a nonvanishing $2\pi$-periodic (resp. $2\pi$-antiperiodic)
 $C^{n-1}$-function.
 Suppose 
$$
\int_{S^1} u(t)\varphi(t)dt=0
\leqno {(3)}
$$
for all $\varphi(t)\in {\cal A}$.
Then the function $u(t)$ changes its sign at least $n+1$
times on $S^1$. }

\medskip
\noindent
{\it Proof.} Suppose that $u$ does not change sign if $n$ is odd and that it
changes
sign only ones if $n$ is even. Let $t_1$ be arbitrary if $n$ is odd and the
zero of
$u$ if $n$ is even. By Corollary A.3 there is a nonvanishing
function $\varphi$ in ${\cal A}$ with a zero in $t_1$ with
 multiplicity $n-1$. We can choose the sign of $\varphi$ such that $
u(t)\varphi(t)\ge0$
for all $t$. Then the integral in (3) being equal to zero implies that
$u\varphi$ vanishes
identically. This is a contradiction since both $u$ and $\varphi$ vanish in
at most one point.
Hence $u$ changes sign at least ones if $n$ is odd and at least twice if $n$
is even.

Now assume that $u$ changes sign only in the
$m$ distinct points
$$
t_1< t_2< \cdots< t_{m}.
$$
where $m< n+1$. Notice that $m$ is even if $n$ is odd and $m$ is
odd if $n$ is even. Hence $n-m$ is an odd integer.
By Corollary A.3, there exists a function
$\psi$ in ${\cal A}$ such that\medskip
\item{(i)} $\psi$ has zero in $t_1$ with multiplicity $n-m$,
\item{(ii)} $\psi$ has a zero in $t_j$ with multiplicity one if $j\ge2$.

\medskip
\noindent
Since the total multiplicity of zeros of
$\psi$ is $n-1$, it follows that $\psi$ has no zeros other than
$t_1,\dots,t_m$.
Since the multiplicities of the zeros at
$t_1,\dots,t_{m}$ are all odd integers, we have that
$\psi$ changes its sign exactly at $t_1,\dots,t_{m}$.
Thus $u\psi$ never changes sign.
The vanishing of the integral in (3), now implies that
$u$ must vanish identically, which is
a contradiction. Hence $u$ changes it sign at least $n+1$ times.
\qed

\medskip
Let ${\cal A}$ be a Chebyshev space of odd (resp.~even) order $n$
in  $C^{n-1}({\bf R})$.
Let $u$ be a nonvanishing $2\pi$-periodic
(resp.~anti $2\pi$-periodic) $C^{n-1}$-function.
Let $s\in S^1$. A function $\varphi_s\in {\cal A}$ is called the
{\it ${\cal A}$-osculating function}
of $u$ at $s$ if the first $n-1$ derivatives
of $\varphi_s$ and $u$ at $s$ coincide.
That $\varphi_s$ exists and is unique follows
from the bijectivity of the mapping $T$ in Theorem A.2 by setting
$t_1=\cdots=t_n=s$.

\medskip
\noindent
{\bf A.5 Definition.}
Suppose ${\cal A}$ and $u$ as above are contained in $C^{n}({\bf R})$.
Then a point $s\in S^1$ is called an {\it  ${\cal A}$-flex}
if the $n$-th derivative of $\varphi_s$ and $u$
coincide in $s$.

\medskip
>From now on we assume that a Chebyshev space ${\cal A}$
of order $n$ is contained in $C^n({\bf
R})$.
We will show  in Theorem A.8
that the number of $\cal A$-flexes on an
$2\pi$-periodic
(resp.~$2\pi$-antiperiodic) function $u$ in $C^n({\bf R})$ is at least
$n+1$. For this we will need the next
proposition which is interesting in
its own right.
It shows that the condition that a Chebyshev space ${\cal A}$ of order $n$
be contained in
$C^{n}({\bf R})$ is a necessary and sufficient condition that it is the
kernel of a disconjugate operator of order $n$.

Notice though that
the concept of a Chebyshev space is essentially
wider than that of kernels of disconjugate
operators. An example of a Chebyshev space of
order $3$ which is in $C^{2}({\bf R})$  but not in
$C^{3}({\bf R})$ can be obtained as follows. Let ${\cal A}_{\gamma}$
be the linear span of the functions $\{1,x,y\}$
where $\gamma=(x,y):S^1\to {\bf R}^2$
is a strictly convex $C^{2}$-curve which is not $C^3$.
Clearly, ${\cal A}_{\gamma}$ is a Chebyshev space of order three since any
line meets the curve $\gamma$
in at most two points counted with multiplicities.

\medskip
\noindent
{\bf A.6 Proposition.}
{\it Let ${\cal A}$ be a Chebyshev space of odd (resp. even) order $n$
in $C^{n+m}({\bf R})$ for $m\ge 0$, then there exists a
unique differential operator
of the form
$$
L_{\cal A}=D^{n}+a_{n-1}D^{n-1}+\cdots +a_{1}D+a_0
$$ 
such that ${\cal A}$ is the kernel of $L_{\cal A}$,
where the coefficients $a_j$ are  $2\pi$-periodic
(resp.~$2\pi$-antiperiodic) $C^m$-functions.}

\medskip
The operator $L_{\cal A}$ in the proposition is called {\it the disconjugate
operator
associated with the Chebyshev space $\cal A$}. The uniqueness of $L_{\cal
A}$ is due to the
fact that its highest order coefficient is normalized to be $1$.

\medskip
\noindent
{\it Proof.}
We fix a point $t\in S^1$ arbitrarily and
define a linear map $T_t:{\cal A}\to {\bf R}^n$ by setting
$$
T_t(\varphi)
=\left(\varphi(t),\varphi'(t),\dots, \varphi^{(n-1)}(t)\right).
$$
By Theorem A.2, the map $T_t$ is bijective.
We define a linear functional $S_t:{\cal A}\to {\bf R}$
by setting $S_t(\varphi)=\varphi^{(n)}(t)$.
Since $S_t\circ T_t^{-1}$ is a linear functional on ${\bf R}^n$,
there exists an $(a_0(t),\dots.,a_{n-1}(t))$ in ${\bf R}^n$ such that
$$
S_t\circ T_t^{-1}(x_0,\dots,x_{n-1})=-\sum_{i=0}^{n-1} a_i(t) x_i,
$$
where the choice of the negative sign in front of the sum will soon become
clear.
It is clear that $(a_0(t),\dots.,a_{n-1}(t))$ is $C^m$
in $t$. Now we have
$$\eqalign{
\varphi^{(n)}(t)&=S_t(\varphi)\cr
&=S_t\circ T_t^{-1}(T_t(\varphi))\cr
&=S_t\circ T_t^{-1}\left(
\varphi(t),\varphi'(t),\dots, \varphi^{(n-1)}(t)\right)\cr
&=-\sum_{i=0}^{n-1} a_i(t) \varphi^{(i)}(t).
}
$$
This proves the existence of the operator $L_{\cal A}$. The uniqueness is
clear.
\qed

The next proposition is a further preparation for the existence of flexes in
Theorem A.8. 

\medskip
\noindent
{\bf A.7 Proposition} {\it
Let ${\cal A}$ be a Chebyshev space of odd (resp.~even) order $n$
in $C^n({\bf R})$. Let $u$ be a
$2\pi$-periodic 
(resp.~$2\pi$-antiperiodic) $C^{n}$-function.
A point $s$ is an ${\cal A}$-flex of $u$ if and only if the
function $L_{\cal A}u$ vanishes at $s$. }

\medskip
\noindent
{\it Proof.} Since $\cal A$ is contained in $C^n({\bf R})$ it has an
associated disconjugate operator $\cal A$ such that $L_{\cal A}\varphi=0$
for all 
$\varphi$ in $\cal A$.
Let $\varphi_s$ in ${\cal A}$ be the osculating function
of $u$ at a point $s$.
Since $L_{\cal A}(\varphi_p)$ vanishes identically,
we have that
$$\eqalign{
u^{(n)}(t)-\varphi_s^{(n)}(t)&=
u^{(n)}(t)+\sum_{i=1}^{n-1} a_i(t) \varphi_s^{(i)}(t) \cr
&=u^{(n)}(t)+\sum_{i=0}^{n-1} a_i(t) u^{(i)}(t) \cr
&=(L_{\cal A}u)(t)
}
$$
for all $t$. Hence $s$ is an $\cal A$-flex if and only if $L_{\cal A}u$
vanishes in
$s$.
\qed

\medskip
\noindent
{\bf A.8 Theorem} {\it
Let ${\cal A}$ be a Chebyshev space of odd (resp.~even) order $n$
in $C^n({\bf R})$. Let $u$ be a
$2\pi$-periodic 
(resp. $2\pi$-antiperiodic) $C^{n}$-function.
Then the number of ${\cal A}$-flexes of $u$ on $S^1$ is at least $n+1$.
}

\medskip
\noindent
{\it Proof.}
The adjoint operator $L^*$ of
a  disconjugate operator $L$ is also disconjugate; see Theorem 9 on p.~104
in [Co].
Let ${\cal A}^*$ be the Chebyshev space of order $n$ corresponding to
$L^*_{\cal A}$ which is also contained in $C^n({\bf R})$.
Then
$$
\int_{S^1} (L_{\cal A}u)(t)\varphi(t)\;dt=\int_{S^1}u(t)(L^*_{\cal
A}\varphi)(t)\;dt=0
$$
for all $\varphi$ in ${\cal A}^*$. We now apply Theorem A.4 to the Chebyshev
space ${\cal A}^*$
and conclude that the function $L_{\cal A}u$ changes sign at least $n+1$
times.
Hence $u$ has at least $n+1$ flexes by Proposition A.7.
\qed

\medskip
Theorem A.8 is optimal; see the example after Theorem 5.1.

\bigskip
\noindent
{\bf Appendix B.}
\medskip

\bigskip
Here we shall prove the following result from Calculus which was
used in the proof of Theorem 5.1.

\medskip
\noindent
{\bf Theorem B.1.} {\it Let $u$ be a $C^n$-function defined on
a neighborhood $I$ of the origin and satisfying
$$ 
u(0)=u'(0)=\cdots=u^{(r)}(0)=0
$$
where $r\le n$.
Then there exists a $C^{n-r-1}$-function
$v$  on $I$
such that
$$
u(t)=t^{r+1} v(t)
$$
for all $t\in I$.
}

\medskip
Applying the following lemma $r+1$ times immediately proves the
theorem.

\medskip
\noindent
{\bf Lemma B.2.} {\it
Let $u$ be a $C^k$-function defined on
a neighborhood $I$ of the origin and satisfying
$u(0)=0$ where $k\ge 1$.
Then there exists a $C^{k-1}$-function $v$ on $I$
such that $u(t)=t\, v(t)$ for all $t\in I$.
}
\medskip
\noindent
{\it Proof.}  We have that
$$
u(t)=\int_0^1 {d u(ts) \over ds} ds
=\int_0^1 t\, u'(ts) ds=t\, v(t)
$$
where 
$$
v(t):=\int_0^1 u'(ts) ds.
$$
It is clear that $v$ is continuous. \qed
\medskip

\noindent{\bf Remark.} Lemma B.2 can be found on p.~89 in [Ar1] with the
redundant assumption
that $u'(0)=0$.

\bigskip\bigskip
\noindent {\bf References}\medskip

\item{[Ar1]} V. Arnold: {\it \'Equations diff\'erentielles ordinaires.\/}
\'Editons Mir, Moscou,
1974.

\item{[Ar2]} V. Arnold: {\it
Remarks on the extatic points of plane curves.
}In: The Gelfand
Mathematical Seminars, 1993--1995, 11--22,
Birkh\"auser, Boston,
1996.

\item{[Ba]} M. Barner:
{\it \"Uber die Mindestanzahl station\"arer Schmiegebenen bei
geschlos\-senen \break strengkonvexen
Raumkurven.}
 Abh.~Math.~Sem.~Univ.~Hamburg {\bf 20} (1956), 196-215.

\item{[Bl]} W. Blaschke: {\it Vorlesungen \"uber Differentialgeometrie II},
Die Grundlehren der Mathematischen Wissenschaften, Band 7, Springer-Verlag,
Berlin 1923. 

\item{[Bo]}
       R. C. Bose:
{\it On the number of circles of curvature perfectly
              enclosing or perfectly enclosed  by
              a closed convex oval.}
 Math. Z.~{\bf 35} (1932),
 16--24.

\item{[Co]} W.A. Coppel: {\it Disconjugacy}. Lecture Notes in Math.~{\bf
220},
 Springer-Verlag, Berlin etc., 1971.

\item{[GMO]}  L. Guieu, E. Mourre \&
V. Yu. Ovsienko: {\it  Theorem on
six vertices of a plane curve via
Sturm theory.} In: The Arnold-Gelfand
Mathematical Seminars, pp.~257--266,
Birkh\"auser, Boston, 1997.

\item{[Ha]}
O. Haupt:
{\it Verallgemeinerung eines Satzes von R.~C. Bose \"uber
              die Anzahl der Schmieg\-kreise eines Ovals, die vom Oval
              umschlossen werden oder das Oval umschlie\ss en.}
J. Reine Angew.~Math.
       {\bf 239/240}
       (1969),
        339--352.

\item{[HK]}
O. Haupt and H. K\"unneth:
{\it Geometrische Ordnungen.}
Die Grundlehren der Mathematischen Wissenschaften, Band 133,
 Springer-Verlag, Berlin, New York,
1967.

\item{[Hy]}
T. Hayashi:
{\it Some geometrical applications of Fourier series.}
 Rend.~Circ.~Mat. Palermo
{\bf  50}
       (1926),
       96--102.

\item{[Ja]}
 S.B. Jackson:
{\it Vertices of plane curves.}
      Bull. Amer. Math. Soc. {\bf 50}
      (1944),
564-578.

\item{[Kn]} H. Kneser: {\it Neuer Beweis
des Vierscheitelsatzes.}
Christiaan Huygens {\bf 2} (1922/23),
315--318.

\item{[Mu1]}
S. Mukhopadhyaya: {\it New methods in the geometry of a plane arc, I.}
Bull. Calcutta Math. Soc. {\bf 1} (1909), 31-- 37. Also in:
Collected geometrical papers, vol. I.   Calcutta University Press,
Calcutta, 1929, 13--20.

\item{[Mu2]} S. Mukhopadhyaya: {\it Sur les
nouvelles m\'ethodes de
g\'eometrie.} C. R. S\'eance Soc.~Math.
France, ann\'ee 1933 (1934),
41-45.

\item{[N\"o]}
 G. N\"obeling: {\it
 \"Uber die Anzahl der ordnungsgeometrischen Scheitel von Kurven II}.
 Geometriae Dedicata  {\bf 31} (1989), 137-149.

\item{[OT]}
V. Ovsienko \& S. Tabachnikov:
{\it Sturm theory, Ghys theorem on zeros of
               the Schwarz\-ian derivative and flattening of
               Legendrian curves}. Selecta Math.~(New Series)
               {\bf 2} (1996), 297-307.
\item{[Ta]}
S. Tabachnikov:
{\it Parametrized plane curves, Minkowski caustics, Minkowski vertices
and conservative line fields}. Enseign.~Math. (2) {\bf 43} (1997), 3--26.

\item{[TU1]} G. Thorbergsson \& M. Umehara: {\it A
unified approach to the
four vertex
theorems, II.} In: Differential  and symplectic  topology of
knots and curves, pp.~229--252, American Math.~Soc.~Transl.~(Series 2) {\bf
190},
Amer.~Math.~Soc., Providence, R.I., 1999.

\item{[TU2]} G. Thorbergsson \& M. Umehara: {\it %
Sextactic points on a simple closed curve.} Preprint.

\item{[Um]} M.
Umehara: {\it A unified approach to the four vertex
theorems, I.}
 In: Differential  and symplectic  topology of
knots and curves, pp.~185--228, American Math.~Soc.~Transl. (Series 2) {\bf
190},
Amer.~Math.~Soc., Providence, R.I., 1999.
\vskip 1.in

\hbox{\parindent=0pt\parskip=0pt
\vbox{\hsize=2.7truein
\obeylines
{
Gudlaugur Thorbergsson
Mathematisches Institut
Universit\"at zu K\"oln
Weyertal 86 - 90
50931 K\"oln
Germany
}\medskip
gthorbergsson@mi.uni-koeln.de
}\hskip 1.5truecm
\vbox{\hsize=3.7truein
\obeylines
{
Masaaki Umehara
Department of Mathematics
Graduate School of Science,
Hiroshima University
Higashi-Hiroshima 739-8526
Japan}\medskip
umehara@math.sci-hiroshima-u.ac.jp
}
}
\bye